\documentclass[english,12pt,onecolumn]{IEEEtran}

\usepackage{babel}
\usepackage{epstopdf}
\usepackage{url}
\usepackage{marvosym}
\usepackage{amssymb}
\usepackage{tikz}
\usepackage{textcomp}
\usepackage{marvosym}
\usepackage[cmex10]{amsmath}
\usepackage{amsmath, bm}
\usepackage{tikz}

\newtheorem{theorem}{Theorem}
\newtheorem{lemma}{Lemma}
\newtheorem{corollary}{Corollary}
\newtheorem{assumption}{Assumption}
\newtheorem{definition}{Definition}
\newtheorem{example}{Example}
\newtheorem{remark}{Remark}

\newcommand{\argmin}{\operatornamewithlimits{argmin}}
\DeclareMathOperator{\id}{Id}
\DeclareMathOperator{\sat}{sat}
\DeclareMathOperator{\diag}{diag}

\begin{document}

\title{Singular Perturbation: When the Perturbation Parameter Becomes a State-Dependent Function\thanks{This work was supported in part by National Key R \& D Program of China (2022ZD0119901), in part by Science and Technology Projects of Liaoning Province (2022JH2/101300238), in part by Liaoning Revitalization Talents Program XLYC2203076, in part by NSFC grants 62325303 and 62333004, and in part by NSF grant ECCS-2210320.}}
\author{Tengfei Liu\thanks{State Key Laboratory of Synthetical Automation for Process Industries, Northeastern University, Shenyang, 110004, China Email: {\tt tfliu@mail.neu.edu.cn}} and Zhong-Ping Jiang\thanks{Department of Electrical and Computer Engineering, New York University, Five MetroTech Center, Brooklyn, NY 11201, USA Email: {\tt zjiang@nyu.edu}}}
\date{\ }

\maketitle

\begin{abstract}
This paper presents a new systematic framework for nonlinear singularly perturbed systems in which state-dependent perturbation functions are used instead of constant perturbation coefficients. Under this framework, general results are obtained for the global robust stability and input-to-state stability of nonlinear singularly perturbed systems. Interestingly, the proposed methodology provides innovative solutions beyond traditional singular perturbation theory for emerging control problems arising from nonlinear integral control, feedback optimization, and formation-based source seeking.
\end{abstract}

\begin{IEEEkeywords}
Singular perturbation, perturbation functions, input-to-state stability (ISS).
\end{IEEEkeywords}

\section{Introduction}
\label{section.introduction}

\subsection{Background and Literature Review}

Many complex systems exhibit multiple time scales. Understanding and utilizing the properties of different time scales and their interactions can significantly simplify system analysis and lead to new systematic design methods \cite{Kokotovic84}. Singular perturbation theory, which uses perturbation coefficients to capture time-scale separation, has been developed to solve multi-time-scale problems related to the integration of information processing and physical objects.

Among the numerous results, significant efforts have been devoted to understanding how fast dynamics can impact the stability properties of slow dynamics; see, for example, \cite{Tikhonov48,Levinson50,Wasow65,Hoppensteadt66,Yoshizawa66,Hahn67,Habets74,Saberi-Khalil84,Kokotovic-Khalil-Oreilly86,OMalley91,Isidori95,Khalil02book,Teel-Moreau-Nesic03,Vecchio-Slotine-TAC-2013}. Results for stochastic systems, time-delay systems, hybrid systems and systems modeled by partial differential equations have also been developed accordingly. See, e.g., \cite{Fridman-Auto-2002,Donchev-Dontchev-SIAMCON-2003,Kabanov-Pergamenshchikov-SIAMCON-1997,Wang-Zhang-Yin-SIAMCON-2006,Yang-Zhang-Sun-Chai-TAC-2011,Tang-Prieur-Girard-TAC-2016,Chen-Yuan-Zheng-TAC-2013}. The theory has found a wide range of applications in areas such as optimal regulation \cite{Shinar-Auto-1983}, high-gain control \cite{Young-Kokotovic-Utkin-TAC-1977,Sannuti-Auto-1983,Marino-IJC-1985}, integral control \cite{Riedle-Kokotovic-TAC-1986}, adaptive control \cite{Ioannou-Kokotovic-Auto-1985}, stochastic control \cite{Kabanov-Pergamenshchikov-SIAMCON-1997}, and extremum seeking \cite{Krstic-Wang-Auto-2000}.

When the stability properties of different time scales are characterized by Lyapunov functions, sum-type or more general composite Lyapunov functions are commonly employed for stability analysis of singularly perturbed systems \cite{Klimushchev-Krasovskii61,Hoppensteadt68,Saberi-Khalil84,Corless-Glielmo92,Corless-Garofalo-Glielmo93,Subotic-Gross-Colombino-Dorfler-TAC-2021}. Alternative methods based on Lyapunov functions include level-set-based analysis in \cite{Habets74,Watbled05}.
Converse Lyapunov theorems \cite{Khalil02book,Sontag-Wang95,Lin-Sontag-Wang96,Clarke-Ledyaev-Stern-JDE-1998,Hahn67} are essential for robustness analysis of nonlinear uncertain systems that lack explicit Lyapunov functions. When it comes to singularly perturbed systems where the stability properties of subsystems are not formulated by explicit Lyapunov functions, converse Lyapunov theorems have been used in \cite{Corless-Glielmo92,Christofides-Teel96}. Additionally, \cite{Moreau-Aeyels-TAC-2000} proposed a Lyapunov-based approach, and \cite{Teel-Moreau-Nesic03} developed a trajectory-based approach for robustness analysis of time-varying singularly perturbed systems with solutions continuously depending on a small parameter.

Among the aforementioned results, the concept of input-to-state stability (ISS) developed by \cite{Sontag89,Sontag-Wang96} and its variants have been utilized to describe the influence of external disturbances and the interaction between the dynamics of different time scales; see, for example, \cite{Christofides-Teel96,Teel-Moreau-Nesic03,Forni-Angeli19}.
It has been proven that if the perturbation coefficient is sufficiently small, the ISS of a singularly perturbed system can be derived from the ISS of the reduced-order subsystem that corresponds to the slow dynamics, possibly in a semi-global practical manner.

We would like to mention that small-gain theorems have been used as a fundamental idea in the analysis of singularly perturbed systems. We recall the insightful discussion by \cite{Christofides-Teel96} on applying the small-gain thinking developed by \cite{Zames66,Mareels-Hill92,Jiang-Teel-Praly94,Teel-Praly95} to stability analysis of singularly perturbed systems. Indeed, the Lyapunov-based conditions in \cite{Saberi-Khalil84} are closely related to small-gain conditions in an $L_2$ setting. As a typical application of the singular perturbation theory, several extremum seeking results, such as \cite{Tan-Nesic-Mareels06,Haring-Johansen17,Hazeleger-Haring-Wouw20,Wang-Qin-Liu-Jiang20}, have been developed with the nonlinear small-gain theorem \cite{Jiang-Teel-Praly94} as a tool for stability analysis. Motivated by emerging control applications such as feedback optimization and formation-based source seeking, we aim to develop a new framework of singular perturbation where we replace the constant perturbation parameter by state-dependent functions.

\subsection{Main Contributions of This Paper}

A singularly perturbed system involving two time-scales often takes the following form:
\begin{align}
\dot{x}&=g_s(x,z,d)\label{plant.trad.slow}\\
\epsilon \dot{z}&=g_f(z,x,w)\label{plant.trad.fast}
\end{align}
where $x\in\mathbb{R}^n$ and $z\in\mathbb{R}^m$ are the states of the subsystems, $d\in\mathbb{R}^p$ and $w\in\mathbb{R}^q$ represent disturbances affecting the subsystems, $g_s:\mathbb{R}^n\times\mathbb{R}^m\times\mathbb{R}^p\rightarrow\mathbb{R}^n$ and $g_f:\mathbb{R}^m\times\mathbb{R}^n\times\mathbb{R}^q\rightarrow\mathbb{R}^m$ describe the system dynamics, and $\epsilon$ is a small positive constant characterizing the different time scales of the two subsystems. By convention, $\epsilon$ is called the perturbation parameter. In the literature (see, e.g., \cite{Kokotovic84}), the model \eqref{plant.trad.slow}--\eqref{plant.trad.fast} is said to be in the standard form if the equilibrium of interest corresponds to an isolated real root of $0=g_f(z,x,0)$, denoted as $z_e = \varphi(x)$.

By setting a new time scale $\tau=t/\delta$ with constant $\delta>0$, the system \eqref{plant.trad.slow}--\eqref{plant.trad.fast} can be rewritten as
\begin{align}
\frac{dx}{d\tau}&=\delta g_s(x,z,d),\label{plant.trad.slow2}\\
\frac{dz}{d\tau}&=\frac{\delta}{\epsilon} g_f(z,x,w).\label{plant.trad.fast2}
\end{align}

However, the standard singular perturbation model \eqref{plant.trad.slow}--\eqref{plant.trad.fast}, or its equivalent counterpart \eqref{plant.trad.slow2}--\eqref{plant.trad.fast2}, exhibit limitations in several emerging applications such as feedback optimization and formation-based source seeking. As shown in Section \ref{section.motivation}, these applications necessitate the tuning of time scales in a nonlinear manner for improved performance, to which the conventional singular perturbation is not directly applicable. With these observations in mind, in this paper, we consider a generalized form of \eqref{plant.trad.slow2}--\eqref{plant.trad.fast2}:
\begin{align}
\dot{x}&=\rho_s(x,z,d)g_s(x,z,d)\label{plant.slow}\\
\dot{z}&=\rho_f(z,x,w)g_f(z,x,w)\label{plant.fast}
\end{align}
which uses positive state-dependent functions $\rho_s:\mathbb{R}^n\times\mathbb{R}^m\times\mathbb{R}^p\rightarrow\mathbb{R}_+$ and $\rho_f:\mathbb{R}^m\times\mathbb{R}^n\times\mathbb{R}^q\rightarrow\mathbb{R}_+$ to describe the time-scale characteristics of the subsystems. For existence and uniqueness of solutions, all the discussions in this paper are based on the assumption that $\rho_s(x,z,d)g_s(x,z,d)$ and $\rho_f(z,x,w)g_f(z,x,w)$ are locally Lipschitz with respect to $(x,z)$ and they are piecewise continuous with respect to $d$ and $w$, respectively. Clearly, if $\rho_s\equiv \delta$ and $\rho_f\equiv \delta/\epsilon$, then the system \eqref{plant.slow}--\eqref{plant.fast} is reduced to \eqref{plant.trad.slow2}--\eqref{plant.trad.fast2}.
With the more general model \eqref{plant.slow}--\eqref{plant.fast}, we expect to improve the design flexibility by tuning the time scales in some nonlinear fashion.

The study presented in this paper adopts the idea of time-scale separation \cite{Tikhonov48} and utilizes the notion of ISS (see \cite{Sontag89}) to characterize the stability properties of the subsystems. Specifically, we assume that the boundary-layer subsystem, i.e., for any fixed $x$, the $z$-subsystem \eqref{plant.fast} with $\rho_f\equiv 1$ is ISS with respect to the equilibrium $z_e$ and the input $w$. It is also assumed that the reduced-order $x$-subsystem \eqref{plant.slow} with $\rho_s\equiv 1$ is ISS with $d$ and $z-z_e$ 
(the error state of the boundary-layer subsystem) as the inputs. We use Lyapunov-based ISS gains (see \cite{Sontag-Wang96}), which depend on perturbation functions, to describe the interaction between the reduced-order and the boundary-layer subsystems, and derive monotonicity conditions on the perturbation functions to guarantee the ISS of the singularly perturbed system. The Lyapunov-based formulation helps simplify the robustness analysis. The proofs of the main results rely on Lyapunov-based ISS small-gain arguments (see \cite{Jiang-Mareels-Wang96}).

The refined singular perturbation results provide useful tools for controller design for complex systems beyond conventional singularly perturbed systems. Applications given in this paper include integral control, feedback optimization, and formation-based source seeking for nonlinear uncertain systems. Moreover, the ISS-based formulation is general in the sense that the main result of this paper can be easily extended to regional and practical stability cases. However, due to space limitations, related research findings will be reported in a companion paper.

\subsection*{Notations and Terminology}

In this paper, we use $|\cdot|$ to denote the Euclidean norm for real vectors and the induced $2$-norm for real matrices. A continuous function $\alpha:\mathbb{R}_+\rightarrow\mathbb{R}_+$ is said to be of class $\mathcal{P}$, denoted by $\alpha\in\mathcal{P}$, if $\alpha(r)>0$ for all $r>0$ and $\alpha(0)\geq 0$; it is said to be of class $\mathcal{PD}$, denoted by $\alpha\in\mathcal{PD}$, if $\alpha(r)>0$ for all $r>0$ and $\alpha(0)=0$.
A continuous function $\alpha:\mathbb{R}_+\rightarrow\mathbb{R}_+$ is said to be of class $\mathcal{K}$, denoted by $\alpha\in\mathcal{K}$, if it is strictly increasing and $\alpha(0)=0$; it is said to of class $\mathcal{K}_\infty$, denoted by $\alpha\in\mathcal{K}_{\infty}$, if it is of class $\mathcal{K}$ and unbounded. A continuous function $\beta:\mathbb{R}_+\times\mathbb{R}_+\rightarrow\mathbb{R}_+$ is said to be of class $\mathcal{KL}$, denoted by $\beta\in\mathcal{KL}$, if, for each fixed $t\in\mathbb{R}_+$, $\beta(\cdot,t)$ is of class $\mathcal{K}$ and, for each fixed $s>0$, $\beta(s,\cdot)$ is a decreasing function and satisfies $\lim_{t\rightarrow\infty}\beta(s,t)=0$.
For $\alpha_1,\alpha_2\in\mathcal{P}$, $\alpha_1\circ\alpha_2$ denotes the composition function, i.e., $\alpha_1\circ\alpha_2(r)=\alpha_1(\alpha_2(r))$ for $r\geq0$. For a function $\varphi:\mathbb{R}^p\rightarrow\mathbb{R}$, $\nabla\varphi$ represents the gradient wherever it exists. $\id$ represents the identity function. $\diag\{a_1,\ldots,a_m\}$ denotes the diagonal matrix with scalars $a_1,\ldots,a_m$ on the diagonal.

\section{Motivational Examples}
\label{section.motivation}

This section gives three examples that motivate the study of this paper.

\subsection{Nonlinear Integral Feedback}
\label{section.motivation.subsection.integralcontrol}

Consider the plant
\begin{align}
\dot{z}^o&=g_f^o(z^o,x^o)\label{integralcontrol.plant1}\\
y^o&=h_f^o(z^o)\label{integralcontrol.plant2}
\end{align}
where $z^o\in\mathbb{R}^m$ is the state, $y^o\in\mathbb{R}^n$ is the output, $x^o\in\mathbb{R}^n$ is the reference input, and $g_f^o:\mathbb{R}^m\times\mathbb{R}^n\rightarrow\mathbb{R}^m$ and $h_f^o:\mathbb{R}^m\rightarrow\mathbb{R}^n$ are locally Lipschitz functions, representing the dynamics and the output map of the plant, respectively.

Suppose that for each fixed $x^o$, the $z^o$-subsystem admits an equilibrium at $\varphi^o(x^o)$ with $\varphi^o:\mathbb{R}^n\rightarrow\mathbb{R}^m$ being a locally Lipschitz function satisfying
\begin{align}
g_f^o(\varphi^o(x^o),x^o)=0\label{integralcontrol.steadystate}
\end{align}
for all $x^o\in\mathbb{R}^n$. Intuitively, $\varphi^o$ represents the steady-state input-state map. Moreover, for each fixed $x^o$, the plant is assumed to be asymptotically stable at the equilibrium $\varphi^o(x^o)$.

Integral negative feedback updates the reference input $x^o$ by integrating the output $y^o$:
\begin{align}
\dot{x}^o=-\epsilon y^o=-\epsilon h_f^o(z^o)\label{integralcontrol.integrator}
\end{align}
where $\epsilon$ is a positive constant. The closed-loop system composed of \eqref{integralcontrol.plant1} and \eqref{integralcontrol.integrator} is in the form of the standard singular perturbation model \eqref{plant.trad.slow}--\eqref{plant.trad.fast}. Based on the conventional singular perturbation theory, in addition to monotonicity conditions on the input-output map $h_f^o\circ\varphi^o$, (local) exponential stability at $\varphi^o(x^o)$ for each fixed $x^o$ is normally required for the plant \eqref{integralcontrol.plant1} to guarantee asymptotic stability of the closed-loop system; see, e.g., \cite{Desoer-Lin-TAC-1985} as well as the recent results \cite{Simpson-Porco21,Lorenzetti-Weiss22}.

The requirement of exponential stability on the plant restricts the applicability of integral negative feedback. Here is an elementary counter-example:
\begin{align}
\dot{z}^o&=g_f^o(z^o,x^o)=-(z^o-x^o)^3\label{integralcontrol.example.plant1}\\
y^o&=h_f^o(z^o)=z^o\label{integralcontrol.example.plant2}\\
\dot{x}^o&=-\epsilon y^o\label{integralcontrol.example.integrator}
\end{align}
Clearly, as a direct application of Lyapunov's first theorem, this system is unstable at the zero equilibrium.

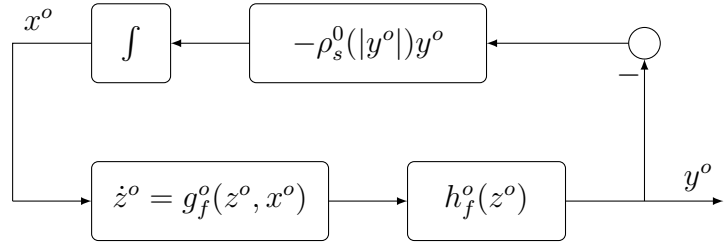
\begin{figure}[h!]
\centering
\begin{tikzpicture}[>=latex,font=\normalsize,scale=2.1]
\draw[->] (0mm,2.5mm) -- (5mm,2.5mm);
\draw [rounded corners=0.9mm] (5mm,0mm) rectangle (20mm,5mm);
\draw (12.5mm,2.5mm) node {$\dot{z}^o=g_f^o(z^o,x^o)$};
\draw[->] (20mm,2.5mm) -- (25mm,2.5mm);
\draw [rounded corners=0.9mm] (25mm,0mm) rectangle (35mm,5mm);
\draw (30mm,2.5mm) node {$h_f^o(z^o)$};
\draw[->] (35mm,2.5mm) -- (45mm,2.5mm);
\draw (45mm,4mm) node [left] {$y^o$};
\draw[->] (40mm,2.5mm) -- (40mm,11.5mm);
\draw (40mm,12.5mm) circle [radius=1mm];
\draw (39mm,10.5mm) node {$-$};
\draw[->] (39mm,12.5mm) -- (30mm,12.5mm);
\draw [rounded corners=0.9mm] (15mm,10mm) rectangle (30mm,15mm);
\draw (22.5mm,12.5mm) node {$-\rho_s^0(|y^o|)y^o$};
\draw[->] (15mm,12.5mm) -- (10mm,12.5mm);
\draw [rounded corners=0.9mm] (5mm,10mm) rectangle (10mm,15mm);
\draw (7.5mm,12.5mm) node {$\int$};
\draw (5mm,12.5mm) -- (0mm,12.5mm) -- (0mm,2.5mm);
\draw (0mm,14mm) node [right] {$x^o$};
\end{tikzpicture}
\caption{Configuration of a control system with nonlinear integral feedback.}
\label{figure.integralcontrol}
\end{figure}

On the basis of the above observations, we propose to use a nonlinear integral control algorithm instead of \eqref{integralcontrol.integrator}:
\begin{align}
\dot{x}^o=-\rho_s^0(|y^o|)y^o=-\rho_s^0(|h_f^o(z^o)|)h_f^o(z^o)\label{integralcontrol.integrator'}
\end{align}
with $\rho_s^0\in\mathcal{P}$. A schematic diagram of a control system with nonlinear integral feedback is shown in Figure \ref{figure.integralcontrol}.

The closed-loop system composed of the plant \eqref{integralcontrol.example.plant1}--\eqref{integralcontrol.example.plant2} and the modified integral controller \eqref{integralcontrol.integrator'} is globally asymptotically stable at the origin, as long as the function $\rho_s^0$ is chosen to be positive definite and satisfy
\begin{align}
\rho_s^0(r)\leq kr^2,~~~\forall r\geq 0,
\end{align}
with $k$ being a positive constant; see Example \ref{example.feedbackoptimization}. Note that the modified integral control system is in the general singular perturbation form \eqref{plant.slow}--\eqref{plant.fast} with state $(x,z)$ coresponding to $(x^o,z^o)$. With the help of a refined singular perturbation theorem, we provide a positive solution to the nonlinear integral control problem; see Section \ref{section.applications.subsection.integralcontrol}.

\subsection{Feedback Optimization}

Feedback optimization is a new direction in control theory that aims to solve an optimization problem with the physical model in the loop, by means of the gradient of the objective function computed in real time. For instance, see \cite{Jokic-Lazar-Bosch2009,Wang-Elia2011,Francois-Bonvin13,Zhao-Topcu-Li-Low14,Brunner-Durr-Ebenbauer2016,Dorfler-Simpson-Porco-Bullo2016,Gan-Low16,Tang-Dvijotham-Low2017,Molzahn17,DallAnese-Simonetto18,Chang-Colombino-Cortes-DallAnese2019,Colombino-DallAnese-Bernstein2020,Hauswirth-Bolognani-Hug-Dorfler2021,Krishnamoorthy-Skogestad22,Liu-Qin-Hong-Jiang22}.

Still consider the plant \eqref{integralcontrol.plant1}--\eqref{integralcontrol.plant2}, which admits a steady-state input-output map $\varphi^o$ defined by \eqref{integralcontrol.steadystate}. Given a continuously differentiable objective function $\Phi:\mathbb{R}^p\times\mathbb{R}^n\rightarrow\mathbb{R}$, under specific convexity conditions, a feedback optimization algorithm updates the reference input $x^o$ such that the state $(x^o,z^o)$ keeps bounded and converges to $(x_*^o,z_*^o)$, a set-point satisfying the steady-state input-output constraint $z_*^o=\varphi^o(x_*^o)$ and minimizing $\Phi(h_f^o(z^o),x^o)$. In such a problem setting,
\begin{align}
x_*^o=\argmin_{x^o}\Phi(h_f^o(\varphi^o(x^o)),x^o).
\end{align}

In \cite{Colombino-DallAnese-Bernstein2020,Hauswirth-Bolognani-Hug-Dorfler2021}, modified variable-metric gradient-flow algorithms \cite{Powell1983,Boyd-Vandenberghe2004} have been used for feedback optimization. Here is an example of such algorithm:
\begin{align}
\dot{x}^o=\rho_s^0(|g_s^o(x^o,z^o)|)g_s^o(x^o,z^o)\label{gradientflow}
\end{align}
where
\begin{align}
g_s^o(x^o,z^o)=-\frac{\partial\Phi(y^o,x^o)}{\partial y^o}\frac{\partial h_f^o(\varphi^o(x^o))}{\partial x^o}-\frac{\partial\Phi(y^o,x^o)}{\partial x^o}
\end{align}
and $\rho_s^0\in\mathcal{P}$ represents the variable metric term. Here, $(x^o,z^o)$ is considered as the arguments of the function $g_s^o$ because $y^o=\varphi^o(z^o)$.

With the method proposed in this paper, we show that the variable metric term indeed helps solve the feedback optimization problem in the absence of the usual exponential stability property for the boundary-layer system at the steady state; see Section \ref{section.applications.subsection.feedbackoptimization}.

\subsection{Formation-Based Source Seeking}
\label{section.motivation.subsection.sourceseeking}

In a typical scenario of source seeking, the objective is to search the position $p_*\in\mathbb{R}^n$ of a signal source based on real-time measurement of the value of the field generated by the signal source. One may model the field with an objective function $h:\mathbb{R}^n\rightarrow\mathbb{R}$ with an extremum point coinciding with the signal source.

Coordination of mobile agents is promising to solve such problems \cite{Zhang-Leonard-TAC-2010}. Consider a group of $N$ mobile agents, with each agent $i$ modeled by
\begin{align}
\dot{p}_i=v_i\label{multiagent}
\end{align}
where $p_i\in\mathbb{R}^n$ is the position and $v_i\in\mathbb{R}^n$ is the velocity (considered as the control input). Given the objective function $h$, each agent $i$ has access to the value $h(p_i)$ of the objective function. With a coordinated control strategy, it is desired that the positions of the agents keep bounded, and the average position
\begin{align}
p_0=\frac{1}{N}\sum_{i=1}^Np_i\label{averageposition}
\end{align}
converges to a neighborhood of the optimal point $p_*$.

The following distributed control law is proposed in this paper to solve the problem:
\begin{align}
v_i=v_i^e+v_i^f\label{coordinatedsourceseekingcontrollaw}
\end{align}
where $v_i^f$ and $v_i^e$ are used for formation control and source seeking, respectively.

In particular, the formation control part of the control law is defined by
\begin{align}
v_i^f=\frac{1}{N}\sum_{j=1}^Na_{ij}(p_i-p_j-d_{ij})\label{formationcontrollaw}
\end{align}
where the $a_{ij}$'s are nonnegative constants, and $d_{ij}$ are constant vectors representing the desired relative position between agents $i$ and $j$ and satisfying $d_{ij}=-d_{ji}$ and
\begin{align}
\det\left(\sum_{j=1}^Nd_{j0}d_{j0}^T\right)>0
\end{align}
with $d_{j0}=\sum_{i=1}^Nd_{ji}/N$ represents the relative position between agent $j$ and the average position in the desired formation. The source-seeking control part of the control law is defined by
\begin{align}
v_i^e=-c_0\sigma\left(\left(\sum_{j=1}^Nd_{j0}d_{j0}^T\right)^{-1}\delta_i\right)\label{formationsourceseekinglaw}
\end{align}
where $c_0$ is a positive constant, $\sigma:\mathbb{R}^n\rightarrow\mathbb{R}^n$ is a saturated function satisfying $\sigma(r)r>0$ for $r\in\mathbb{R}^n\backslash\{0\}$ and $\sigma(0)=0$, and $\delta_i$ is generated by the following distributed averaging algorithm
\begin{align}
\dot{\delta}_i&=-(\delta_i-Nd_{i0}h(p_i))-\sum_{j=1}^Nb_{ij}(q_i-q_j)\label{gradientestimation1}\\
\dot{q}_i&=\mu\sum_{j=1}^Nb_{ij}(\delta_i-\delta_j)\label{gradientestimation2}
\end{align}
with $\mu$ being a positive constant and $b_{ij}$ being nonnegative constants.

\begin{figure}[h!]
\centering
\begin{tikzpicture}[>=latex,font=\normalsize,scale=2.1]
\draw[white,fill=gray,opacity=0.15] (0mm,0mm) circle (15mm);
\draw[white,fill=gray,opacity=0.15] (2.8mm,-0.3mm) circle (11mm);
\draw[white,fill=gray,opacity=0.15] (4.5mm,-0.5mm) circle (8mm);
\draw[white,fill=gray,opacity=0.15] (5.5mm,-0.7mm) circle (6mm);
\draw[white,fill=gray,opacity=0.15] (6.2mm,-0.9mm) circle (4mm);
\draw[white,fill=gray,opacity=0.15] (7mm,-1.7mm) circle (2mm);
\draw (7mm,-1.7mm) node {$\diamond$};
\draw (7mm,-1.7mm) node [right] {$p_*$};
\draw[densely dotted] (-12mm,5mm) -- ++(-70:13mm) -- ++(-190:13mm) -- cycle;
\draw[->] (-12mm,5mm) -- ++(-3:9mm);
\draw (-12mm,5mm) ++(-3:9mm) ++(-3:-0.5mm) node [right] {$v_i$};
\draw[densely dashed,->] (-12mm,5mm) -- ++(-93:3.276mm);
\draw (-12mm,5mm) ++(-93:3.276mm) ++(-93:1.5mm) node {$v_i^f$};
\draw[densely dashed,->] (-12mm,5mm) -- ++(17:9.578mm);
\draw (-12mm,5mm)++(17:9.578mm) ++(17:1mm) node {$v_i^e$};
\draw[densely dashed] (-12mm,5mm) ++(-93:3.276mm) -- ++(17:9.578mm);
\draw[densely dashed] (-12mm,5mm) ++(17:9.578mm) -- ++(-93:3.276mm);
\draw[->] (-12mm,5mm) ++(-70:13mm) -- ++(-1:9mm);
\draw (-12mm,5mm) ++(-70:13mm) ++(-1:9mm) ++(-1:-0.5mm) node [right] {$v_j$};
\draw[->] (-12mm,5mm) ++(-130:13mm) -- ++(2:9mm);
\draw (-12mm,5mm) ++(-130:13mm) ++(2:9mm) ++(2:-0.5mm) node [right] {$v_k$};
\draw[fill=white,thick] (-12mm,5mm) circle (1mm);
\draw[fill=white,thick] (-12mm,5mm) ++(-70:13mm) circle (1mm);
\draw[fill=white,thick] (-12mm,5mm) ++(-130:13mm) circle (1mm);
\draw (-13mm,5mm) node [left] {$\dot{p}_i=v_i$};
\end{tikzpicture}
\caption{Formation-Based Source Seeking.}
\label{figure.formationbasedsourceseeking}
\end{figure}
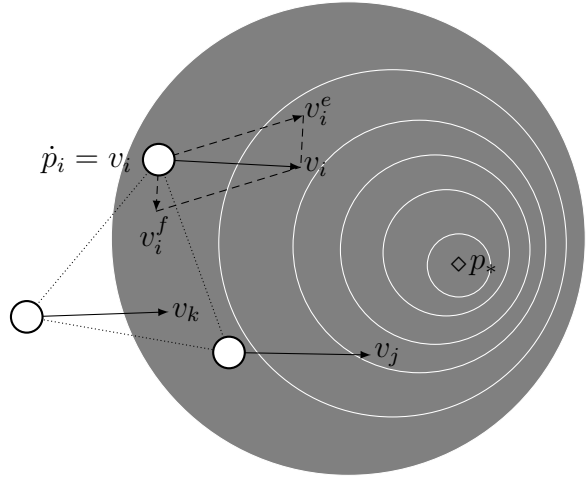

The basic idea of the formation-based source-seeking algorithm is shown in Figure \ref{figure.formationbasedsourceseeking}.

Indeed, the averaging algorithm gives an estimate of the gradient information of the objective function. In the ideal case of fixed $h(p_i)$ for $i=1,\ldots,N$, the averaging algorithm results in $\lim_{t\rightarrow\infty}\delta_j(t)=\sum_{i=1}^Nd_{i0}h(p_i)$ for $j=1,\ldots,N$ with $\sum_{i=1}^Nd_{i0}h(p_i)$ being a mean-value-based estimate of the gradient value $\nabla h(p_0)$. This is the basic idea of employing the averaging algorithm. As far as we know, distributed online-averaging algorithms have not been utilized for gradient-free source-seeking in past literature. It can be recognized that the control law only uses the relative positions between the agents and the measured value of the objective function, and does not rely on the global information of $p_*$ or the gradient measurement of the objective function.

The processes of gradient estimation and source seeking interact with each other in a nonlinear fashion. Intuitively, the gradient estimation process should be fast enough. But it is nontrivial to find an explicit condition to guarantee the convergence of the closed-loop states. By considering the two processes as fast and slow subsystems, respectively, we propose a novel systematic design based on the refined singular perturbation theorem; see Section \ref{section.applications.subsection.sourceseeking}.

\section{Problem Formulation}
\label{section.problemformulation}

Figure \ref{figure.singularperturbationequivalence} shows the idea of considering a singularly perturbed system as an interconnection of the reduced-order and the boundary-layer subsystems.

\begin{figure}[h!]
\centering
\begin{tikzpicture}[>=latex,font=\normalsize,scale=2.1]
\draw [rounded corners=0.9mm] (5mm,0mm) rectangle (15mm,5mm);
\draw (10mm,2.5mm) node {$z$};
\draw[->] (-5mm,1.5mm) -- (5mm,1.5mm);
\draw[->] (0mm,3.5mm) -- (5mm,3.5mm);
\draw (15mm,2.5mm) -- (20mm,2.5mm) -- (20mm,11.5mm);
\draw[->] (20mm,11.5mm) -- (15mm,11.5mm);
\draw[->] (25mm,13.5mm) -- (15mm,13.5mm);
\draw [rounded corners=0.9mm] (5mm,10mm) rectangle (15mm,15mm);
\draw (10mm,12.5mm) node {$x$};
\draw (5mm,12.5mm) -- (0mm,12.5mm) -- (0mm,3.5mm);
\draw (25mm,15mm) node [left] {$d$};
\draw (-5mm,0mm) node [right] {$w$};
\draw[white] (-5mm,-12.5mm) rectangle (25mm,27.5mm);
\end{tikzpicture}
\begin{tikzpicture}[>=latex,scale=2.1]
\draw (0mm,15mm) node {$\Leftrightarrow$};
\draw[white] (-2mm,-5mm) rectangle (2mm,35mm);
\end{tikzpicture}
\begin{tikzpicture}[>=latex,font=\normalsize,scale=2.1]
\draw [rounded corners=0.9mm] (5mm,0mm) rectangle (15mm,5mm);
\draw (10mm,2.5mm) node {$z$};
\draw[->] (-5mm,1.5mm) -- (5mm,1.5mm);
\draw[->] (0mm,3.5mm) -- (5mm,3.5mm);
\draw (15mm,2.5mm) -- (20mm,2.5mm);
\draw[->] (20mm,26.5mm) -- (15mm,26.5mm);
\draw[->] (25mm,28.5mm) -- (15mm,28.5mm);
\draw [rounded corners=0.9mm] (5mm,25mm) rectangle (15mm,30mm);
\draw (10mm,27.5mm) node {$x$};
\draw (5mm,27.5mm) -- (0mm,27.5mm) -- (0mm,3.5mm);
\draw (25mm,30mm) node [left] {$d$};
\draw (-5mm,0mm) node [right] {$w$};
\draw[->] (0mm,10mm) -- (6.5mm,10mm);
\draw [rounded corners=0.9mm] (6.5mm,8mm) rectangle (13.5mm,12mm);
\draw (10mm,10mm) node {$\varphi(x)$};
\draw[->] (13.5mm,10mm) -- (19.3mm,10mm);
\draw[->] (0mm,20mm) -- (6.5mm,20mm);
\draw [rounded corners=0.9mm] (6.5mm,18mm) rectangle (13.5mm,22mm);
\draw (10mm,20mm) node {$\varphi(x)$};
\draw[->] (13.5mm,20mm) -- (19.3mm,20mm);
\draw[->] (20mm,2.5mm) -- (20mm,9.3mm);
\draw (20mm,10mm) circle [radius=0.7mm];
\draw (19mm,8mm) node {$+$};
\draw (18mm,9mm) node {$-$};
\draw[->] (20mm,10.7mm) -- (20mm,19.3mm);
\draw (20mm,20mm) circle [radius=0.7mm];
\draw (19mm,18mm) node {$+$};
\draw (18mm,19mm) node {$+$};
\draw (20mm,20.7mm) -- (20mm,26.5mm);
\draw[densely dashed] (-5mm,15mm) -- (25mm,15mm);
\draw[white] (-5mm,-5mm) rectangle (25mm,35mm);
\draw (10mm,-3mm) node {Boundary-Layer Subsystem};
\draw (10mm,33mm) node {Reduced-Order Subsystem};
\end{tikzpicture}
\caption{Singularly perturbed system \eqref{plant.slow}--\eqref{plant.fast} as an interconnected system, where $\varphi$ represents the steady-state map of the $z$-subsystem.}
\label{figure.singularperturbationequivalence}
\end{figure}

Following the convention, the following assumptions are made on the steady state and the stability properties of the subsystems.

Assumption \ref{assumption.steadystate} means that $(x,z)=(0,\varphi(0))$ is an equilibrium of the system \eqref{plant.slow}--\eqref{plant.fast} with $d\equiv 0$ and $w\equiv 0$.

\begin{assumption}\label{assumption.steadystate}
There exists a locally Lipschitz function $\varphi:\mathbb{R}^n\rightarrow\mathbb{R}^m$ such that
\begin{align}
g_s(0,\varphi(0),0)&=0\label{assumption.steadystate.p1}\\
g_f(\varphi(x),x,0)&=0\label{assumption.steadystate.p2}
\end{align}
for all $x\in\mathbb{R}^n$.
\end{assumption}

Assumption \ref{assumption.stability} gives the stability properties of the subsystems corresponding to different time scales.

\begin{assumption}\label{assumption.stability}
There exist continuously differentiable functions $V_s:\mathbb{R}^n\rightarrow\mathbb{R}_+$ and $V_f:\mathbb{R}^m\times\mathbb{R}^n\rightarrow\mathbb{R}_+$ such that
\begin{itemize}
\item for the reduced-order subsystem \eqref{plant.slow} with $\rho_s\equiv 1$, there exist $\underline{\alpha}_s,\overline{\alpha}_s,\gamma_s,\chi_s\in\mathcal{K}_{\infty}$, and $\alpha_s\in\mathcal{PD}$ such that
\begin{align}
&\underline{\alpha}_s(|x|)\leq V_s(x)\leq\overline{\alpha}_s(|x|)\label{assumption.slowsubsystem.Lya1}\\
&V_s(x)\geq\max\left\{\gamma_s(V_f(z,x)),\chi_s(|d|)\right\}\Rightarrow\frac{\partial V_s(x)}{\partial x}g_s(x,z,d)\leq-\alpha_s(V_s(x))\label{assumption.slowsubsystem.Lya2}
\end{align}
for all $z\in\mathbb{R}^m$, $x\in\mathbb{R}^n$ and $d\in\mathbb{R}^p$;
\item for the boundary-layer subsystem \eqref{plant.fast} with $\rho_f\equiv 1$, there exist $\underline{\alpha}_f,\overline{\alpha}_f,\chi_f\in\mathcal{K}_{\infty}$, and $\alpha_f,\lambda_{f1},\lambda_{f2}\in\mathcal{PD}$ such that
\begin{align}
&\underline{\alpha}_f(|z-\varphi(x)|)\leq V_f(z,x)\leq\overline{\alpha}_f(|z-\varphi(x)|)\label{assumption.fastsubsystem.Lya1}\\
&V_f(z,x)\geq\chi_f(|w|)\Rightarrow\frac{\partial V_f(z,x)}{\partial z}g_f(z,x,w)\leq-\alpha_f(V_f(z,x))\label{assumption.fastsubsystem.Lya2}\\
&\left|\frac{\partial V_f(z,x)}{\partial x}\right|\leq\lambda_{f1}(V_f(z,x))+\lambda_{f2}(V_s(x))\label{assumption.fastsubsystem.Lya3}
\end{align}
for all $z\in\mathbb{R}^m$, $x\in\mathbb{R}^n$ and $w\in\mathbb{R}^q$.
\end{itemize}
\end{assumption}

\begin{remark}
Assumption \ref{assumption.stability} is a Lyapunov characterization of the ISS properties of the unperturbed reduced-order and the unperturbed boundary-layer subsystems, i.e., with the perturbation functions identically equal to one. The cross terms between the subsystems are directly described by Lyapunov functions instead of state variables, which simplifies the stability analysis of the singularly perturbed system.
\end{remark}

\begin{remark}
In the absence of disturbances $d$ and $w$, it is usual to assume an asymptotically stable (AS) equilibrium for the reduced-order subsystem \eqref{plant.slow} with $z=\varphi(x)$ and $\rho_s\equiv 1$, in accordance with Assumption \ref{assumption.stability}. For the more general case in which the systems are subject to exogenous disturbances, the ISS property of the reduced-order subsystem with $V_f(z,x)$ as the input (or equivalently, with $z-\varphi(x)$ as the input) demanded by Assumption \ref{assumption.stability} is stronger. However, when we limit ourselves to the local case, the corresponding ISS assumption is equivalent to the $0$-AS assumption (i.e., AS of the zero-input system at the origin) due to the equivalence between $0$-AS and local ISS (see \cite{Sontag-Wang96}). In addition, it should be mentioned that the ISS assumption is weaker than assumptions made in the past literature, such as the exponential stability and additional regularity conditions required by \cite{Corless-Glielmo92}, as ISS does not guarantee exponential convergence.
\end{remark}

We aim to derive readily-checkable conditions on the perturbation functions $\rho_s$ and $\rho_f$ to guarantee the ISS of the generalized singularly perturbed system \eqref{plant.slow}--\eqref{plant.fast}.

\section{Main Results}
\label{section.main}

Theorem \ref{theorem.main} presents a condition on $\rho_s$ and $\rho_f$ to ensure ISS of the system \eqref{plant.slow}--\eqref{plant.fast}.

\begin{theorem}\label{theorem.main}
Under Assumptions \ref{assumption.steadystate} and \ref{assumption.stability}, the system \eqref{plant.slow}--\eqref{plant.fast} is ISS with $(x,z-\varphi(0))$ as the state and $(d,w)$ as the input, if $\rho_s$ and $\rho_f$ satisfy the following conditions:
\begin{enumerate}
\item there exists a $\underline{\rho}_s\in\mathcal{P}$ such that
\begin{align}
&V_s(x)\geq\max\left\{\gamma_s(V_f(z,x)),\chi_s(|d|)\right\}\Rightarrow\rho_s(x,z,d)\geq\underline{\rho}_s(V_s(x))\label{theorem.main.con1}
\end{align}
for all $x\in\mathbb{R}^n$, $z\in\mathbb{R}^m$ and $d\in\mathbb{R}^p$;
\item there exist $\gamma_f\in\mathcal{K}_{\infty}$ and $\overline{\rho}_s,\underline{\rho}_f\in\mathcal{P}$ satisfying
\begin{align}
\gamma_f\circ\gamma_s(r)&<r\label{theorem.main.con2}\\
\overline{\rho}_s(r)(\lambda_{f1}(r)+\lambda_{f2}\circ\gamma_f^{-1}(r))&<\underline{\rho}_f(r)\alpha_f(r)\label{theorem.main.con3}
\end{align}
for all $r>0$, such that
\begin{align}
&V_f(z,x)\geq\max\left\{\gamma_f(V_s(x)),\gamma_f\circ\chi_s(|d|),\chi_f(|w|)\right\}\nonumber\\
&\Rightarrow\begin{cases}
\rho_s(x,z,d)|g_s(x,z,d)|\leq\overline{\rho}_s(V_f(z,x)),\\
\rho_f(z,x,w)\geq\underline{\rho}_f(V_f(z,x))
\end{cases}\label{theorem.main.con4}
\end{align}
for all $x\in\mathbb{R}^n$, $z\in\mathbb{R}^m$, $d\in\mathbb{R}^p$ and $w\in\mathbb{R}^q$.
\end{enumerate}
Moreover, for certain $\sigma\in\mathcal{K}_{\infty}$ being continuously differentiable on $(0,\infty)$ and satisfying $\gamma_s(r)<\sigma(r)<\gamma_f^{-1}(r)$ for all $r>0$, the following function
\begin{align}
V(x,z)=\max\left\{V_s(x),\sigma(V_f(z,x))\right\}\label{theorem.main.Lyapunovconstruction}
\end{align}
is positive definite and radially unbounded with respect to $(x,z-\varphi(x))$, and satisfies
\begin{align}
&V(x,z)\geq\max\left\{\chi_s(|d|),\sigma\circ\chi_f(|w|)\right\}\nonumber\\
&\Rightarrow\frac{\partial V(x,z)}{\partial x}\rho_s(x,z,d)g_s(x,z,d)+\frac{\partial V(x,z)}{\partial z}\rho_f(z,x,w)g_f(z,x,w)\leq-\alpha(V(x,z))\label{theorem.main.Lyapunovconstructionderivative}
\end{align}
wherever $V$ is differentiable, with $\alpha\in\mathcal{PD}$.
\end{theorem}

\begin{IEEEproof}
We treat the singularly perturbed system as an interconnection of the reduced-order and the boundary-layer subsystems. We first show that each subsystem features specific ISS properties when associated with a perturbation function, and then perform a small-gain analysis leading to the ISS of the interconnected system.

A straightforward implication of property \eqref{assumption.slowsubsystem.Lya2} in Assumption \ref{assumption.stability} and condition \eqref{theorem.main.con1} is that
\begin{align}
&V_s(x)\geq\max\left\{\gamma_s(V_f(z,x)),\chi_s(|d|)\right\}\Rightarrow\frac{\partial V_s(x)}{\partial x}\rho_s(x,z,d)g_s(x,z,d)\leq-\alpha_{bs}(V_s(x))\label{slowgain}
\end{align}
for all $x\in\mathbb{R}^n$, $z\in\mathbb{R}^m$ and $d\in\mathbb{R}^p$, with
\begin{align}
\alpha_{bs}(r)=\underline{\rho}_s(r)\alpha_s(r)
\end{align}
for $r\geq 0$. Since $\underline{\rho}_s\in\mathcal{P}$ and $\alpha_s\in\mathcal{PD}$, we have $\alpha_{bs}\in\mathcal{PD}$.

In the case of
\begin{align}
V_f(z,x)\geq\max\left\{\gamma_f(V_s(x)),\gamma_f\circ\chi_s(|d|),\chi_f(|w|)\right\},
\end{align}
we have
\begin{align}
&\frac{\partial V_f(z,x)}{\partial z}\rho_f(z,x,w)g_f(z,x,w)+\frac{\partial V_f(z,x)}{\partial x}\rho_s(x,z,d)g_s(x,z,d)\nonumber\\
&\leq-\rho_f(z,x,w)\alpha_f(V_f(z,x))+\lambda_{f1}(V_f(z,x))\rho_s(x,z,d)|g_s(x,z,d)|+\lambda_{f2}(V_s(x))\rho_s(x,z,d)|g_s(x,z,d)|\nonumber\\
&\leq-\underline{\rho}_f(V_f(z,x))\alpha_f(V_f(z,x))+\lambda_{f1}(V_f(z,x))\overline{\rho}_s(V_f(z,x))+\lambda_{f2}(\gamma_f^{-1}(V_f(z,x)))\overline{\rho}_s(V_f(z,x))\nonumber\\
&=:-\alpha_{bf}(V_f(z,x))\label{p0}
\end{align}
for all $x\in\mathbb{R}^n$, $z\in\mathbb{R}^m$, $d\in\mathbb{R}^p$ and $w\in\mathbb{R}^q$, where we have used properties \eqref{assumption.fastsubsystem.Lya2} and \eqref{assumption.fastsubsystem.Lya3} in Assumption \ref{assumption.stability} for the first inequality, and used condition \eqref{theorem.main.con4} for the second inequality. That is,
\begin{align}
&V_f(z,x)\geq\max\left\{\gamma_f(V_s(x)),\gamma_f\circ\chi_s(|d|),\chi_f(|w|)\right\}\Rightarrow\nonumber\\
&\frac{\partial V_f(z,x)}{\partial z}\rho_f(z,x,w)g_f(z,x,w)+\frac{\partial V_f(z,x)}{\partial x}\rho_s(x,z,d)g_s(x,z,d)\leq-\alpha_{bf}(V_f(z,x))\label{fastgain}
\end{align}
for all $x\in\mathbb{R}^n$, $z\in\mathbb{R}^m$, $d\in\mathbb{R}^p$ and $w\in\mathbb{R}^q$. Condition \eqref{theorem.main.con3} together with $\alpha_f,\overline{\rho}_s\in\mathcal{PD}$ guarantees that $\alpha_{bf}\in\mathcal{PD}$.

Properties \eqref{slowgain} and \eqref{fastgain} describe the Lyapunov-based ISS properties of the reduced-order subsystem with $x$ as the state and the boundary-layer subsystem with $z-\varphi(x)$ as the state, respectively, and characterize their interconnection by ISS gains $\gamma_f$ and $\gamma_s$. Condition \eqref{theorem.main.con2} ensures the satisfaction of the nonlinear small-gain condition \cite{Jiang-Mareels-Wang96}. By applying the Lyapunov-based ISS small-gain theorem (refer to Theorem \ref{theorem.Lyapunovsmallgain} in the Appendix), we can prove the ISS property of the system \eqref{plant.slow}--\eqref{plant.fast} with $(x,z-\varphi(x))$, or equivalently $(x,z-\varphi(0))$, as the state and $(d,w)$ as the input, and construct an ISS-Lyapunov function in \eqref{theorem.main.Lyapunovconstruction} with the Lyapunov-ISS property \eqref{theorem.main.Lyapunovconstructionderivative}. This ends the proof of Theorem \ref{theorem.main}.
\end{IEEEproof}

\begin{remark}
Given Assumptions \ref{assumption.steadystate} and \ref{assumption.stability}, it is always possible to find perturbation functions that guarantee ISS of the singularly perturbed system. Specifically, the condition of Theorem \ref{theorem.main} can be satisfied by setting $\rho_f\equiv 1$ and $\underline{\rho}_f(\cdot)\equiv 1$, selecting $\underline{\rho}_s, \overline{\rho}_s\in\mathcal{PD}$ small enough, and choosing $\rho_s$ solely depending on $|x|$ and small enough. Alternatively, by setting $\rho_s\equiv 1$ and selecting $\rho_f$ depending on $|z-\varphi(x)|$ and large enough, the condition can also be satisfied. However, in practice, a more detailed analysis may be necessary when the perturbation functions are restricted to some predetermined forms. See Section \ref{section.special} below.
\end{remark}

\begin{remark}
Using small-gain techniques in singular perturbation analysis is not a new idea. In \cite{Christofides-Teel96}, small-gain arguments \cite{Zames66, Mareels-Hill92, Jiang-Teel-Praly94, Teel-Praly95} are employed to analyze the stability of singularly perturbed systems. The Lyapunov-based conditions in \cite{Saberi-Khalil84} are closely related to small-gain conditions in an $L_2$ framework. Nonlinear small-gain theorems have been widely used in applications of the singular perturbation theory, such as in extremum seeking studies, as seen in \cite{Tan-Nesic-Mareels06, Haring-Johansen17, Hazeleger-Haring-Wouw20, Wang-Qin-Liu-Jiang20}. It is expected that the methods proposed in this paper pave the foundation to further generalize the above-mentioned results.
\end{remark}

\begin{remark}
In the boundedness and convergence analysis in \cite{Tikhonov48,Hoppensteadt66}, as well as subsequent studies such as \cite{Watbled05}, level sets are defined by coordinately evaluating the behaviors of the Lyapunov functions of the reduced-order and boundary-layer subsystems, which is somewhat related to the idea of max-type Lyapunov functions. Additionally, in \cite{Habets74}, a Lyapunov-like function is constructed by taking the maximum of linearly weighted Lyapunov-like functions of the two subsystems with different time scales; see \cite[Theorems 6 and 7]{Habets74}. To handle the state-dependent perturbation functions, this paper constructs a Lyapunov function by taking the maximum of nonlinearly weighted Lyapunov functions of the subsystems corresponding to different time scales.
\end{remark}

For certain cases, constant-valued $\rho_s$ and $\rho_f$ exist for ISS of the generalized singularly perturbed system \eqref{plant.slow}--\eqref{plant.fast}. Before presenting the next main theorem for such cases, we define two functions, $\breve{g}_s$ and $\bar{g}_s$. Given $g_s$ and $\varphi$ satisfying Assumption \ref{assumption.steadystate}, there exists a $\breve{g}_s\in\mathcal{K}_{\infty}$ such that
\begin{align}
|g_s(x,z,d)|\leq \breve{g}_s(|x|+|z-\varphi(x)|+|d|)\label{gupperbound00}
\end{align}
for all $x\in\mathbb{R}^n$, $z\in\mathbb{R}^m$ and $d\in\mathbb{R}^q$. Then, define $\bar{g}_s:\mathbb{R}_+\rightarrow\mathbb{R}_+$ as
\begin{align}
\bar{g}_s(r)=\breve{g}_s(\underline{\alpha}_s^{-1}\circ\gamma_s(r)+\underline{\alpha}_f^{-1}(r)+\gamma_s\circ\chi_s^{-1}(r))\label{gupperbound0}
\end{align}
for $r\geq 0$. Clearly, $\bar{g}_s\in\mathcal{K}_{\infty}$.

The following theorem gives a sufficient condition for the existence of constant-valued $\rho_s$ and $\rho_f$ for ISS of the generalized singularly perturbed system.

\begin{theorem}\label{theorem.mainconstant}
Under Assumptions \ref{assumption.steadystate} and \ref{assumption.stability}, suppose that there exists a positive constant $c_0>0$ satisfying
\begin{align}
c_0\bar{g}_s(r)\bar{\lambda}_f(r)<\alpha_f(r)\label{theorem.mainconstant.condition}
\end{align}
for all $r>0$, with $\bar{g}_s$ defined by \eqref{gupperbound0} and $\bar{\lambda}_f=\lambda_{f1}+\lambda_{f2}\circ\gamma_s$, and set $\rho_s\equiv c_s$ and $\rho_f\equiv c_f$ with positive constants $c_s$ and $c_f$ satisfying 
\begin{align}
\frac{c_s}{c_f}=c_0.
\end{align}
Then, the system \eqref{plant.slow}--\eqref{plant.fast} is ISS with $(x,z-\varphi(0))$ as the state and $(d,w)$ as the input.
\end{theorem}

\begin{IEEEproof}
Theorem \ref{theorem.mainconstant} is proved by verifying the satisfaction of the conditions of Theorem \ref{theorem.main}.

We first consider the case with $\rho_s\equiv c_0$ and $\rho_f\equiv 1$. Set $\underline{\rho}_s\equiv c_0$ and $\underline{\rho}_f\equiv 1$. This guarantees the satisfaction of condition \eqref{theorem.main.con1} and the second implication in condition \eqref{theorem.main.con4} in Theorem \ref{theorem.main}.

Condition \eqref{theorem.mainconstant.condition} guarantees the existence of $\hat{\gamma}_s\in\mathcal{K}_{\infty}$ satisfying $\hat{\gamma}_s>\gamma_s$ such that
\begin{align}
c_0\hat{g}_s(r)\hat{\lambda}_f(r)<\alpha_f(r)
\end{align}
holds with
\begin{align}
\hat{g}_s(r)&=\breve{g}_s(\underline{\alpha}_s^{-1}\circ\hat{\gamma}_s(r)+\underline{\alpha}_f^{-1}(r)+\hat{\gamma}_s\circ\chi_s^{-1}(r))\label{theorem.mainconstant.proof.hatgs}\\
\hat{\lambda}_f(r)&=\lambda_{f1}(r)+\lambda_{f2}\circ\hat{\gamma}_s(r)
\end{align}
for $r\geq 0$.

Choose
\begin{align}
\gamma_f=\hat{\gamma}_s^{-1},~~~~~~
\overline{\rho}_s(\cdot)=c_0\hat{g}_s(\cdot).
\end{align}
Then, conditions \eqref{theorem.main.con2} and \eqref{theorem.main.con3} are satisfied.

Recalling the definition of $\breve{g}_s$ in \eqref{gupperbound00} and using properties \eqref{assumption.slowsubsystem.Lya1} and \eqref{assumption.fastsubsystem.Lya1} in Assumption \ref{assumption.stability}, we have
\begin{align}
|g_s(x,z,d)|\leq \breve{g}_s(\underline{\alpha}_s^{-1}(V_s(x))+\underline{\alpha}_f^{-1}(V_f(z,x))+|d|)
\end{align}
for all $x\in\mathbb{R}^n$, $z\in\mathbb{R}^m$ and $d\in\mathbb{R}^p$. Moreover, with $\rho_s\equiv c_0$, in the case of
\begin{align}
V_f(z,x)\geq\max\left\{\gamma_f(V_s(x)),\gamma_f\circ\chi_s(|d|)\right\},
\end{align}
it holds that
\begin{align}
&\rho_s(x,z,d)|g_s(x,z,d)|\nonumber\\
&\leq c_0\breve{g}_s(\underline{\alpha}_s^{-1}\circ\gamma_f^{-1}(V_f(z,x))+\underline{\alpha}_f^{-1}(V_f(z,x))+\gamma_f^{-1}\circ\chi_s^{-1}(V_f(z,x)))\nonumber\\
&=c_0\hat{g}_s(V_f(z,x))=\overline{\rho}_s(V_f(z,x)),
\end{align}
where we have used $\gamma_f=\hat{\gamma}_s^{-1}$ for the inequality and used the definition of $\hat{g}_s$ in \eqref{theorem.mainconstant.proof.hatgs} for the equality. Thus, the satisfaction of the first implication in condition \eqref{theorem.main.con4} is verified. This proves the case of $\rho_s\equiv c_0$ and $\rho_f\equiv 1$ of Theorem \ref{theorem.mainconstant}.

Note that linearly adjusting the time scale does not change the stability property of a singularly perturbed system in the standard form. This means that the singularly perturbed system is still ISS by setting $\rho_s\equiv c_s$ and $\rho_f\equiv c_f$ as long as $c_s/c_f=c_0$. This ends the proof of Theorem \ref{theorem.mainconstant}.
\end{IEEEproof}

We use the following elementary example to illustrate the effectiveness of Theorem \ref{theorem.mainconstant}.

\begin{example}\label{example.saturateddynamics}
Consider a nonlinear system with saturated dynamics:
\begin{align}
\dot{x}&=-c_0\sat(z)\label{example.saturateddynamics.plant1}\\
\dot{z}&=-\sat(z-x)\label{example.saturateddynamics.plant2}
\end{align}
where $x\in\mathbb{R}$ and $z\in\mathbb{R}$ are the states, $c_0$ is a positive constant, and $\sat$ is the saturation function defined by $\sat(r)=\min\{1,\max\{-1,r\}\}$ for $r\in\mathbb{R}$.

Clearly, the system \eqref{example.saturateddynamics.plant1}--\eqref{example.saturateddynamics.plant2} is in the form of \eqref{plant.slow}--\eqref{plant.fast} with
\begin{align}
g_s(x,z)&=-\sat(z),~~~g_f(z,x)=-\sat(z-x).
\end{align}
The constant $c_0$ is considered as the perturbation coefficient. It appears to be nontrival to formulate the stability properties of the subsystems with Lyapunov functions featuring quadratic supply functions.

Direct calculation verifies the satisfaction of Assumption \ref{assumption.steadystate} with
\begin{align}
\varphi(r)=r.
\end{align}
By defining Lyapunov functions
\begin{align}
V_s(x)=\frac{1}{2}x^2,~~~V_f(z,x)=\frac{1}{2}(z-x)^2,\label{example.saturateddynamics.Lyapunovfunctions}
\end{align}
we can verify the satisfaction of Assumption \ref{assumption.stability} with
\begin{align}
&\underline{\alpha}_s(r)=\overline{\alpha}_s(r)=\underline{\alpha}_f(r)=\overline{\alpha}_f(r)=\frac{1}{2}r^2,\\
&\gamma_s(r)=\frac{1}{\o^2}r,~~~\alpha_s(r)=\frac{1}{1-\o}\psi\circ\overline{\alpha}_s^{-1}(r),\\
&\alpha_f(r)=\psi\circ\overline{\alpha}_f^{-1}(r),~~~\lambda_{f1}(r)=\sqrt{2r},~~~\lambda_{f2}(r)=0,
\end{align}
for $r\geq 0$, where $\psi(r)=r\sat(r)$, and $\o$ is a constant that can be chosen arbitrarily from the interval $(0,1)$.

In accordance with \eqref{gupperbound00} and \eqref{gupperbound0}, define
\begin{align}
\breve{g}_s(r)=\sat(r),~~~\bar{g}_s(r)=\sat\left(\sqrt{\frac{2}{\o^2}r}+\sqrt{2r}\right).
\end{align}

Then, condition \eqref{theorem.mainconstant.condition} in Theorem \ref{theorem.mainconstant} is equivalent to
\begin{align}
c_0\sat\left(\sqrt{\frac{2}{\o^2}r}+\sqrt{2r}\right)\sqrt{2r}<\sqrt{2r}\sat\left(\sqrt{2r}\right)
\end{align}
i.e.,
\begin{align}
c_0\sat\left(\frac{1+\o}{\o}\sqrt{2r}\right)<\sat\left(\sqrt{2r}\right)\label{example.saturateddynamics.condition0}
\end{align}
for all $r>0$. Note that $c_0\sat(r)\leq\sat(c_0r)$ for all $0<c_0<1$ and $r\geq 0$. Thus, condition \eqref{example.saturateddynamics.condition0} is satisfied if
\begin{align}
0<c_0<\frac{\o}{1+\o}.\label{example.saturateddynamics.condition}
\end{align}
One can find a constant $\o$ from the interval $(0,1)$ to satisfy condition \eqref{example.saturateddynamics.condition} as long as
\begin{align}
0<c_0<\frac{1}{2}.\label{example.saturateddynamics.condition'}
\end{align}
With Theorem \ref{theorem.mainconstant}, the singularly perturbed system is globally asymptotically stable (GAS) at the origin if $c_0$ satisfies \eqref{example.saturateddynamics.condition'}.
\end{example}

Corollary \ref{corollary.perturbationcoefficient} follows from Theorem \ref{theorem.mainconstant}, for the case where the comparison functions in Assumption \ref{assumption.stability} are compatible with globally Lipschitz dynamics and quadratic Lyapunov functions.

\begin{corollary}\label{corollary.perturbationcoefficient}
Under Assumptions \ref{assumption.steadystate} and \ref{assumption.stability}, there exist constant-valued $\rho_s,\rho_f\in\mathcal{P}$ such that the system \eqref{plant.slow}--\eqref{plant.fast} is ISS with $(x,z-\varphi(0))$ as the state and $(d,w)$ as the input, if
\begin{itemize}
\item $g_s$, $\gamma_s$, $\alpha_f^{-1}$ are globally Lipschitz functions,
\item $\underline{\alpha}_f$ and $\underline{\alpha}_s$ are square functions, and
\item $\lambda_{f1}$ and $\lambda_{f2}$ are square root functions.
\end{itemize}
\end{corollary}

\begin{IEEEproof}
Consider the following two cases. If $\chi_s$ is not less than a square function, then, with the conditions of the corollary satisfied, $\bar{g}_s$ defined by \eqref{gupperbound0} is not larger than some square root function. The global Lipschitz continuity of $\alpha_f^{-1}$ guarantees that $\alpha_f$ is lower bounded by a strictly increasing linear function. Recall that $\lambda_f$ is a square root function. Then, condition \eqref{theorem.mainconstant.condition} can be satisfied by choosing $c_0$ small enough.

If $\chi_s$ is not large enough, then we choose $\chi_s'>\chi_s$ such that $\chi_s'$ is not less than a square function. In this case, Assumption \ref{assumption.stability} still holds with $\chi_s$ replaced by $\chi_s'$. Then, we can still perform the analysis for the first case. This ends the proof of Corollary \ref{corollary.perturbationcoefficient}.
\end{IEEEproof}

\begin{remark}
Corollary \ref{corollary.perturbationcoefficient} is consistent with the results from \cite{Corless-Glielmo92} obtained under exponential stability assumptions. In particular, the ISS condition in Assumption \ref{assumption.stability} can be relaxed to exponential stability at the origin due to the robustness of $0$-exponentially stable systems with globally Lipschitz dynamics \cite[Lemma 4.6]{Khalil02book}. Moreover, if the conditions of Corollary \ref{corollary.perturbationcoefficient} are satisfied, in particular, if the gain function $\gamma_s\in\mathcal{K}_{\infty}$ is globally Lipschitz, then a constant $\sigma_0$ can be found such that $\sigma(r)=\sigma_0r$ defined for $r\geq 0$ satisfies the condition of Theorem \ref{theorem.main}. Consequently, the Lyapunov function $V$ defined in \eqref{theorem.main.Lyapunovconstruction} takes the form
\begin{align}
V(x,z)=\max\{V_s(x),\sigma_0V_f(z,x)\}.
\end{align}
This means that a Lyapunov function can be constructed by taking the maximum of linearly weighted Lyapunov functions of the subsystems, which is consistent with the so-called consistency analysis in \cite{Habets74}.
\end{remark}

\section{Perturbation Functions in Special Forms}
\label{section.special}

This section focuses on specific cases where one of the perturbation functions takes a constant value, while the other one depends on the magnitude of the dynamics of its corresponding subsystem.

\subsection{Tuning the Time Scale of the Reduced-Order Subsystem}

For the singularly perturbed system \eqref{plant.slow}--\eqref{plant.fast}, we consider
\begin{align}
\rho_s(x,z,d)&=\rho_s^0(|g_s(x,z,d)|)\label{special1}\\
\rho_f(z,x,w)&\equiv 1\label{special2}
\end{align}
with $\rho_s^0\in\mathcal{P}$ to be chosen. Systems in various applications, such as feedback optimization using a variable-metric gradient-flow algorithm in Subsection \ref{section.applications.subsection.feedbackoptimization}, can be expressed in this form.

The following theorem proposes an additional monotonicity condition for the existence of $\rho_s^0$ to guarantee ISS of the special class of singularly perturbed systems.

\begin{theorem}\label{theorem.tuningslow}
Under Assumptions \ref{assumption.steadystate} and \ref{assumption.stability}, if there exists a nondecreasing $\overline{\rho}_s\in\mathcal{P}$ satisfying \eqref{theorem.main.con3} with $\underline{\rho}_f(\cdot)\equiv 1$, then one can find a $\rho_s^0\in\mathcal{P}$ for \eqref{special1} such that the conditions \eqref{theorem.main.con1}, \eqref{theorem.main.con2} and \eqref{theorem.main.con4} are satisfied, and the system \eqref{plant.slow}--\eqref{plant.fast} with $\rho_s$ and $\rho_f$ defined by \eqref{special1}--\eqref{special2} is ISS with $(x,z-\varphi(0))$ as the state and $(d,w)$ as the input.
\end{theorem}

\begin{IEEEproof}
For any $\gamma_s\in\mathcal{K}_{\infty}$, there exists a $\gamma_f\in\mathcal{K}_{\infty}$ satisfying \eqref{theorem.main.con2}.
Given $\gamma_f\in\mathcal{K}_{\infty}$ and nondecreasing $\overline{\rho}_s\in\mathcal{P}$, we now find a $\rho_s^0\in\mathcal{P}$ for the satisfaction of \eqref{theorem.main.con4}.

Recall the definition of $\breve{g}_s\in\mathcal{K}_{\infty}$ in \eqref{gupperbound00}. Then, with Assumption \ref{assumption.stability} satisfied, in the case of
\begin{align}
V_f(z,x)\geq\max\left\{\gamma_f(V_s(x)),\gamma_f\circ\chi_s(|d|),\chi_f(|w|)\right\},\label{theorem.tuningslow.case}
\end{align}
it holds that
\begin{align}
|g_s(x,z,d)|&\leq\breve{g}_s(\underline{\alpha}_s^{-1}(V_s(x))+\underline{\alpha}_f^{-1}(V_f(z,x))+|d|)\nonumber\\
&\leq\hat{g}_s(\underline{\alpha}_s^{-1}(\gamma_f^{-1}(V_f(z,x)))+\underline{\alpha}_f^{-1}(V_f(z,x))+\gamma_f^{-1}(\chi_s^{-1}(V_f(z,x))))\nonumber\\
&=\tilde{g}_s(V_f(z,x))\label{gupperbound2}
\end{align}
where
\begin{align}
\tilde{g}_s(r)=\breve{g}_s(\underline{\alpha}_s^{-1}\circ\gamma_f^{-1}(r)+\underline{\alpha}_f^{-1}(r)+\gamma_f^{-1}\circ\chi_s^{-1}(r))\label{gupperbound}
\end{align}
for $r\geq 0$. Clearly, $\tilde{g}_s\in\mathcal{K}_{\infty}$. 

Given $\overline{\rho}_s\in\mathcal{P}$, we choose $\rho_s^0\in\mathcal{P}$ such that
\begin{align}
\rho_s^0(r)r\leq\overline{\rho}_s(\tilde{g}_s^{-1}(r))\label{rhos0condition}
\end{align}
holds for all $r\geq 0$. Since $\overline{\rho}_s$ and $\tilde{g}_s$ are nondecreasing, in the case of \eqref{theorem.tuningslow.case}, we have
\begin{align}
\rho_s^0(|g_s(x,z,d)|)|g_s(x,z,d)|&\leq\overline{\rho}_s(\tilde{g}_s^{-1}(|g_s(x,z,d)|))\nonumber\\
&\leq\overline{\rho}_s(V_f(z,x)).
\end{align}
This guarantees the satisfaction of \eqref{theorem.main.con4}.

With the $\rho_s^0\in\mathcal{P}$ chosen above, we now find a $\underline{\rho}_s\in\mathcal{P}$ for the satisfaction of \eqref{theorem.main.con1}.

Under Assumptions \ref{assumption.steadystate} and \ref{assumption.stability}, the set
\begin{align}
\Omega(r)=\biggl\{(x,z,d):~|x|=r,~V_s(x)\geq\max\{\gamma_s(V_f(z,x)),\chi_s(|d|)\}\biggr\}
\end{align}
is nonempty and bounded for any $r\geq 0$. Then, we can define
\begin{align}
g_s^l(r)&=\inf\{g_s(x,z,d):(x,z,d)\in\Omega(r)\}\\
g_s^u(r)&=\sup\{g_s(x,z,d):(x,z,d)\in\Omega(r)\}
\end{align}
for $r\geq 0$.

With the satisfaction of Assumptions \ref{assumption.steadystate} and \ref{assumption.stability}, $|x|=0$ and $V_s(x)\geq\max\{\gamma_s(V_f(z,x)),\chi_s(|d|)\}$ together imply $x=0$, $z=\varphi(0)$, $d=0$ and thus $g_s(x,z,d)=g_s(0,\varphi(0),0)=0$. This means $g_s^l(0)=0$. Using Assumption \ref{assumption.stability}, if $|x|=r>0$ and $V_s(x)\geq\max\{\gamma_s(V_f(z,x)),\chi_s(|d|)\}$, then $|g_s(x,z,d)|>0$. Otherwise, property \eqref{assumption.slowsubsystem.Lya2} does not hold. This means $g_s^l(r)>0$ for all $r>0$. Similarly, one can verify $g_s^u(0)=0$, and prove $g_s^u(r)>0$ for all $r>0$ by using $g_s^u(r)\geq g_s^l(r)$ for all $r\geq 0$. Given positive definite $g_s^l$ and $g_s^u$, one can find $g_s^{lc},g_s^{uc}\in\mathcal{PD}$ satisfying $g_s^{lc}(r)\leq g_s^l(r)$ and $g_s^{uc}(r)\geq g_s^u(r)$ for all $r\geq 0$.

Based on the discussion above, it is obvious that
\begin{align}
&V_s(x)\geq\max\{\gamma_s(V_f(z,x)),\chi_s(|d|)\}\nonumber\\
&\Rightarrow g_s^{lc}(|x|)\leq|g_s(x,z,d)|\leq g_s^{uc}(|x|)\nonumber\\
&\Rightarrow\rho_s^0(|g_s(x,z,d)|)\geq\underline{\rho}_s^0(|x|)
\end{align}
with
\begin{align}
\underline{\rho}_s^0(r)=\min\left\{\rho_s^0(\xi):g_s^{lc}(r)\leq\xi\leq g_s^{uc}(r)\right\}
\end{align}
for all $r\geq 0$. Then, condition \eqref{theorem.main.con1} is satisfied by defining
\begin{align}
\underline{\rho}_s(r)=\min\left\{\underline{\rho}_s^0(\xi):\overline{\alpha}_s^{-1}(r)\leq\xi\leq\underline{\alpha}_s^{-1}(r)\right\}
\end{align}
for all $r\geq 0$. Then, Theorem \ref{theorem.tuningslow} can be proved with Theorem \ref{theorem.main}.
\end{IEEEproof}

\begin{remark}
Given $\lambda_f,\alpha_f\in\mathcal{PD}$, it may not always be possible to find a non-decreasing $\overline{\rho}_s\in\mathcal{P}$ to satisfy \eqref{theorem.main.con3} with $\underline{\rho}_f(\cdot)\equiv 1$. For example, consider the case where $\alpha_f(r)=r$ and $\lambda_f(r)=r^2$. However, one can definitely find a non-decreasing $\overline{\rho}_s\in\mathcal{P}$ to locally satisfy \eqref{theorem.main.con3} with $\underline{\rho}_f(\cdot)\equiv 1$. This means that local ISS can always be guaranteed with perturbation functions in the form of \eqref{special1}--\eqref{special2} under the assumptions made in Section \ref{section.problemformulation}.
\end{remark}

\subsection{Tuning the Time Scale of the Boundary-Layer Subsystem}

We suppose
\begin{align}
\rho_s(x,z,d)&\equiv 1\label{special1'}\\
\rho_f(z,x,w)&=\rho_f^0(|g_f(z,x,w)|)\label{special2'}
\end{align}
with $\rho_f^0\in\mathcal{P}$ to be chosen.

The following theorem proposes a lower bound condition on the fast dynamics for the existence of $\rho_f^0$ to guarantee the ISS property of the special class of singularly perturbed systems.

\begin{theorem}\label{theorem.tuningfast}
Under Assumptions \ref{assumption.steadystate} and \ref{assumption.stability}, if there exists a $g_f^l\in\mathcal{K}_{\infty}$ such that
\begin{align}
V_f(z,x)\geq\chi_f(|w|)\Rightarrow|g_f(z,x,w)|\geq g_f^l(|z-\varphi(x)|)\label{theorem.tuningfast.con}
\end{align}
for all $x\in\mathbb{R}^n$, $z\in\mathbb{R}^m$ and $w\in\mathbb{R}^q$, then one can find a $\rho_f^0\in\mathcal{P}$ for \eqref{special2'} such that the conditions \eqref{theorem.main.con1}, \eqref{theorem.main.con2}, \eqref{theorem.main.con3} and \eqref{theorem.main.con4} are satisfied, and the system \eqref{plant.slow}--\eqref{plant.fast} with $\rho_s$ and $\rho_f$ defined by \eqref{special1'}--\eqref{special2'} is ISS with $(x,z-\varphi(0))$ as the state and $(d,w)$ as the input.
\end{theorem}

\begin{IEEEproof}
Condition \eqref{theorem.main.con1} is readily satisfied with $\rho_s\equiv 1$ and $\underline{\rho}_s\equiv 1$. For any $\gamma_s\in\mathcal{K}_{\infty}$, there exists a $\gamma_f\in\mathcal{K}_{\infty}$ satisfying \eqref{theorem.main.con2}.

With property \eqref{gupperbound2} in the proof of Theorem \ref{theorem.tuningslow}, we can ensure the satisfaction of the first implication in \eqref{theorem.main.con4} by setting $\overline{\rho}_s=\tilde{g}_s$ with $\tilde{g}_s$ defined in \eqref{gupperbound2}.

With $\lambda_{f1}$, $\lambda_{f2}$ and $\alpha_f$ given by Assumption \ref{assumption.stability} and $\overline{\rho}_s=\tilde{g}_s$ chosen above, one can always find a nondecreasing $\underline{\rho}_f\in\mathcal{P}$ to satisfy \eqref{theorem.main.con3}.

The second implication in \eqref{theorem.main.con4} is satisfied by choosing a nondecreasing $\rho_f^0\in\mathcal{P}$ such that
\begin{align}
\rho_f^0(g_f^l(r))\geq\underline{\rho}_f(\overline{\alpha}_f(r))\label{theorem.tuningfast.proof.0}
\end{align}
for all $r\geq 0$. Indeed, with condition \eqref{theorem.tuningfast.con} satisfied, $V_f(z,x)\geq\chi_f(|w|)$ implies that
\begin{align}
\rho_f^0(|g_f(z,x,w)|)&\geq\rho_f^0(g_f^l(|z-\varphi(x)|))\nonumber\\
&\geq\underline{\rho}_f(\overline{\alpha}_f(|z-\varphi(x)|))\nonumber\\
&\geq\underline{\rho}_f(V_f(z,x)).\label{theorem.tuningfast.proof.1}
\end{align}
Then, Theorem \ref{theorem.tuningfast} can be proved with Theorem \ref{theorem.main}.
\end{IEEEproof}

\begin{remark}
Assumption \ref{assumption.stability} guarantees the existence of a $g_f^l\in\mathcal{PD}$ for the satisfaction of condition \eqref{theorem.tuningfast.con}. This can be verified by referring to the discussion on the existence of $g_s^l$ for $g_s$ in the proof of Theorem \ref{theorem.tuningslow}. The requirement of $g_f^l\in\mathcal{K}_{\infty}$ is stronger. It means that the boundary-layer subsystem should be fast enough to be further tuned to guarantee the ISS of the singularly perturbed system. Instead of relying on condition \eqref{theorem.tuningfast.con}, one may assume the satisfaction of condition \eqref{theorem.main.con3} of Theorem \ref{theorem.main} with $\overline{\rho}_s=\tilde{g}_s$ and a bounded $\underline{\rho}_f\in\mathcal{PD}$. Under the alternative assumption, following the reasoning in \eqref{theorem.tuningfast.proof.1}, one can still find a non-decreasing $\rho_f^0\in\mathcal{P}$ to satisfy \eqref{theorem.tuningfast.proof.0}.
\end{remark}

\section{Applications}
\label{section.applications}

This section shows that the refined singular perturbation results serve as a tool to address stability issues arising from the three applications briefly discussed in Section \ref{section.motivation}: nonlinear integral control, feedback optimization, and formation-based source seeking.

\subsection{Nonlinear Integral Control}
\label{section.applications.subsection.integralcontrol}

The main result facilitates solving the integral control problem for nonlinear uncertain systems. To simplify the discussion and emphasize our contribution, we consider the disturbance-free plant studied in Section \ref{section.motivation.subsection.integralcontrol} and aim for a GAS result.

In particular, consider the plant \eqref{integralcontrol.plant1}--\eqref{integralcontrol.plant2}, and suppose that for each fixed $x^o$, the $z^o$-subsystem admits an equilibrium at $\varphi^o(x^o)$ with $\varphi^o:\mathbb{R}^n\rightarrow\mathbb{R}^m$ being a locally Lipschitz function satisfying \eqref{integralcontrol.steadystate}. With $x^o$ considered as the reference input, integral control aims to steer the output $y^o$ to zero by appropriately updating $x^o$.

The following assumption is made on the steady-state input-output map $h_f^o\circ\varphi^o$.

\begin{assumption}\label{assumption.integralcontrol.steadystate}
There exists an $x_e^o$ satisfying $h_f^o(\varphi^o(x_e^o))=0$, and there exist $\phi_1\in\mathcal{K}_{\infty}$ and $\phi_2\in\mathcal{PD}$ such that
\begin{align}
&|x^o-x_e^o|\geq\phi_1(|\delta|)\Rightarrow(x^o-x_e^o)^Th_f^o(\varphi^o(x^o)+\delta)\geq\phi_2(|x^o-x_e^o|)\label{integralcontrol.signcondition}
\end{align}
for all $x^o\in\mathbb{R}^n$ and $\delta\in\mathbb{R}^m$.
\end{assumption}

Condition \eqref{integralcontrol.signcondition} means that the steady-state input-output map $h_f^o\circ\varphi^o$ keeps inside the half-space determined by the error state $x^o-x_e^o$, and also guarantees that $x_e^o$ is the unique zero point of $h_f^o\circ\varphi^o$. Condition \eqref{integralcontrol.signcondition} is naturally satisfied if $h_f^o$ is globally Lipschitz and $h_f^o\circ\varphi^o$ satisfies the strong monotonicity condition:
\begin{align}
(x^o-x_e^o)^Th_f^o(\varphi^o(x^o))\geq c|x^o-x_e^o|^2
\end{align}
for all $x^o\in\mathbb{R}^n$, with consant $c>0$. The strong monotonocity condition is indispensible for quite a few integral control results including \cite{Desoer-Lin-TAC-1985} as well as the recent results \cite{Simpson-Porco21,Lorenzetti-Weiss22}.

In addition, we assume GAS of the plant at steady states.

\begin{assumption}\label{assumption.integralcontrol.GAS}
There exists a continuously differentiable Lyapunov function $V_f^o:\mathbb{R}^m\times\mathbb{R}^n\rightarrow\mathbb{R}_+$ such that for each fixed $x^o$, it holds that
\begin{align}
&\underline{\alpha}_f^o(|z^o-\varphi^o(x^o)|)\leq V_f^o(z^o,x^o)\leq\overline{\alpha}_f^o(|z^o-\varphi^o(x^o)|)\label{integralcontrol.fastsubsystem.Lya1}\\
&\frac{\partial V_f^o(z^o,x^o)}{\partial z^o}g_f^o(z^o,x^o)\leq-\alpha_f^o(V_f^o(z^o,x^o))\label{integralcontrol.fastsubsystem.Lya2}\\
&\left|\frac{\partial V_f^o(z^o,x^o)}{\partial x^o}\right|\leq\lambda_{f1}^o(|z^o-\varphi^o(x^o)|)+\lambda_{f2}^o(|x^o-x_e^o|)\label{integralcontrol.fastsubsystem.Lya3}
\end{align}
for all $z^o\in\mathbb{R}^m$, where $\underline{\alpha}_f^o,\overline{\alpha}_f^o\in\mathcal{K}_{\infty}$ and $\alpha_f^o,\lambda_{f1}^o,\lambda_{f2}^o\in\mathcal{PD}$.
\end{assumption}

Under Assumptions \ref{assumption.integralcontrol.steadystate} and \ref{assumption.integralcontrol.GAS}, $(x^o,z^o)=(x_e^o,\varphi^o(x_e^o))$ is the unique equilibrium of the integral control system composed of \eqref{integralcontrol.plant1}--\eqref{integralcontrol.plant2} and \eqref{integralcontrol.integrator'}.

The following theorem shows the validity of the nonlinear integral control law \eqref{integralcontrol.integrator'}.

\begin{theorem}\label{theorem.integralcontrol}
Under Assumptions \ref{assumption.integralcontrol.steadystate} and \ref{assumption.integralcontrol.GAS}, there exists a $\rho_s^0\in\mathcal{P}$ such that the integral control system composed of the plant \eqref{integralcontrol.plant1}--\eqref{integralcontrol.plant2} and the integral control law \eqref{integralcontrol.integrator'} is GAS at the equilibrium $(x^o,z^o)=(x_e^o,\varphi^o(x_e^o))$.
\end{theorem}

\begin{IEEEproof}
We first introduce a coordinate change:
\begin{align}
x=x^o-x_e^o,~~z=z^o,~~y=y^o\label{integralcontrol.coordinatetransformation}
\end{align}
to transform the closed-loop system into
\begin{align}
\dot{x}&=-\rho_s^0(|h_f(z)|)h_f(z)\\
\dot{z}&=g_f(z,x)
\end{align}
where $h_f(z)=h_f^o(z)=h_f^o(z^o)$ and $g_f(z,x)=g_f^o(z,x+x_e^o)=g_f^o(z^o,x^o)$. With $\varphi(x)=\varphi^o(x+x_e^o)$, it can be checked that
\begin{align}
g_s(0,\varphi(0))&=-h_f(\varphi(0))=-h_f^o(\varphi^o(x_e^o))=0\\
g_f(\varphi(x),x)&=g_f^o(\varphi^o(x+x_e^o),x+x_e^o)=0
\end{align}
for all $x\in\mathbb{R}^n$. This verifies the satisfaction of Assumption \ref{assumption.steadystate}.

For the $z$-subsystem, with the coordinate transformation \eqref{integralcontrol.coordinatetransformation}, properties \eqref{integralcontrol.fastsubsystem.Lya1}--\eqref{integralcontrol.fastsubsystem.Lya3} are equivalent to
\begin{align}
&\underline{\alpha}_f(|z-\varphi(x)|)\leq V_f(z,x)\leq\overline{\alpha}_f(|z-\varphi(x)|)\label{integralcontrol.fastsubsystem.Lya1'}\\
&\frac{\partial V_f(z,x)}{\partial z}g_f(z,x)\leq-\alpha_f(V_f(z,x))\label{integralcontrol.fastsubsystem.Lya2'}\\
&\left|\frac{\partial V_f(z,x)}{\partial x}\right|\leq\lambda_{f1}^o(|z-\varphi(x)|)+\lambda_{f2}^o(|x|)\label{integralcontrol.fastsubsystem.Lya3'}
\end{align}
where $V_f(z,x)=V_f^o(z,x+x_e^o)$, $\underline{\alpha}_f=\underline{\alpha}_f^o$, and $\overline{\alpha}_f=\overline{\alpha}_f^o$.

For the $x$-subsystem, we consider a quadratic Lyapunov function candidate
\begin{align}
V_s(x)=\frac{1}{2}x^Tx.
\end{align}
Then, condition \eqref{assumption.slowsubsystem.Lya1} is satisfied with $\underline{\alpha}_s(r)=\overline{\alpha}_s(r)=r^2/2$ for $r\geq 0$. Note that
\begin{align}
\frac{\partial V_s(x)}{\partial x}g_s(x,z)=-x^Th_f(z).
\end{align}
Thus,
\begin{align}
&V_s(x)\geq\frac{1}{2}\phi_1^2(|z-\varphi(x)|)\Rightarrow|x|\geq\phi_1(|z-\varphi(x)|)\nonumber\\
&\Rightarrow\frac{\partial V_s(x)}{\partial x}g_s(x,z)\leq-\phi_2(|x|)\leq-\phi_2\circ\overline{\alpha}_s^{-1}(V_s(x))\label{integralcontrol.slowsubsystem.Lya}
\end{align}
for all $x\in\mathbb{R}^n$ and $z\in\mathbb{R}^m$, where condition \eqref{integralcontrol.signcondition} is used for the last implication.

Define $\lambda_{f1}=\lambda_{f1}^o\circ\underline{\alpha}_f^{-1}$ and $\lambda_{f2}=\lambda_{f2}^o\circ\underline{\alpha}_s^{-1}$, and replace property \eqref{integralcontrol.fastsubsystem.Lya3'} with
\begin{align}
\left|\frac{\partial V_f(z,x)}{\partial x}\right|\leq\lambda_{f1}(V_f(z,x))+\lambda_{f2}(V_s(x)).
\end{align}
Define $\alpha_s=\phi_2\circ\overline{\alpha}_s^{-1}$ and replace property \eqref{integralcontrol.slowsubsystem.Lya} with
\begin{align}
V_s(x)\geq\gamma_s(V_f(z,x))\Rightarrow\frac{\partial V_s(x)}{\partial x}g_s(x,z)\leq-\alpha_s(V_s(x)).
\end{align}
Then, the conditions in Assumption \ref{assumption.stability} are verified.

With the satisfaction of Assumptions \ref{assumption.steadystate} and \ref{assumption.stability}, Theorem \ref{theorem.tuningslow} can be readily used to fine tune $\rho_s^0\in\mathcal{P}$ for GAS of the integral control system at the equilibrium.
\end{IEEEproof}

One may replace $y^o$ in the integral control law \eqref{integralcontrol.integrator'} with the output-tracking error $y^o-y^r$ when the output $y^o$ is expected to be steered to some nonzero point $y^r$.

\subsection{Feedback Optimization}
\label{section.applications.subsection.feedbackoptimization}

Consider the feedback optimization system composed of the plant \eqref{integralcontrol.plant1}--\eqref{integralcontrol.plant2} and the feedback optimization algorithm \eqref{gradientflow}. In this problem setting, Assumptions \ref{assumption.integralcontrol.steadystate} and \ref{assumption.integralcontrol.GAS} are still made on the plant.

Denote
\begin{align}
\breve{\Phi}(x^o)=\Phi(h_f^o(\varphi^o(x^o)),x^o).
\end{align}
The following assumption is made on the objective function.

\begin{assumption}\label{assumption.feedbackoptimization.objectivefunction}
There exists a positive constant $\omega$ such that
\begin{align}
(\xi_1-\xi_2)^T(\nabla\breve{\Phi}(\xi_1)-\nabla\breve{\Phi}(\xi_2))\geq\omega|\xi_1-\xi_2|^2\label{assumption.feedbackoptimization.objectivefunction.prop1}
\end{align}
holds for all $\xi_1,\xi_2\in\mathbb{R}^n$.
\end{assumption}

In addition, it is assumed that $g_s^o$ admits the following Lipschitz property.

\begin{assumption}\label{assumption.feedbackoptimization.Lipschitz}
There exists a positive constant $L_{gs}$ such that
\begin{align}
|g_s^o(x^o,z^o+\tilde{z}^o)-g_s^o(x^o,z^o)|\leq L_{gs}|\tilde{z}^o|
\end{align}
holds for all $x^o\in\mathbb{R}^n$, $z^o\in\mathbb{R}^m$ and $\tilde{z}^o\in\mathbb{R}^m$.
\end{assumption}

Condition \eqref{assumption.feedbackoptimization.objectivefunction.prop1} ensures the existence and uniqueness of a global minimum point $x_*^o$ for $\breve{\Phi}$ (see, e.g., \cite{Nesterov2004}). Then, it can be verified that $(x^o,z^o)=(x_*^o,\varphi^o(x_*^o))$ is the unique equilibrium of the feedback optimization system. Based on the assumptions above, the following theorem gives our result on feedback optimization.

\begin{theorem}\label{theorem.feedbackoptimization}
Under Assumptions \ref{assumption.integralcontrol.steadystate}, \ref{assumption.integralcontrol.GAS}, \ref{assumption.feedbackoptimization.objectivefunction} and \ref{assumption.feedbackoptimization.Lipschitz}, there exists a $\rho_s^0\in\mathcal{P}$ such that the feedback optimization system composed of the plant \eqref{integralcontrol.plant1}--\eqref{integralcontrol.plant2} and the variable-metric gradient-flow algorithm \eqref{gradientflow} is GAS at the equilibrium $(x^o,z^o)=(x_*^o,\varphi^o(x_*^o))$.
\end{theorem}

\begin{IEEEproof}[Sketch of Proof]
Some technical derivations are similar with those for the proof of Theorem \ref{theorem.integralcontrol}, and are omitted here to save space.

We first introduce a coordinate transformation:
\begin{align}
x=x^o-x_*^o,~~z=z^o,~~y=y^o\label{feedbackoptimization.coordinatetransformation}
\end{align}
to further transform the closed-loop system into
\begin{align}
\dot{x}&=\rho_s^0(|g_s(x,z)|)g_s(x,z)\\
\dot{z}&=g_f(z,x)
\end{align}
with $g_s(x,z)=g_s^o(x+x_*^o,z)=g_s^o(x^o,z^o)$ and $g_f(z,x)=g_f^o(z,x+x_*^o)=g_f^o(z^o,x^o)$. Define $\varphi(x)=\varphi^o(x^o)=\varphi^o(x+x_*^o)$. Then, it can be checked that
\begin{align}
g_s(0,\varphi(0))&=g_s^o(x_*^o,\varphi^o(x_*^o))=0\\
g_f(\varphi(x),x)&=g_f^o(\varphi^o(x+x_*^o),x+x_*^o)=0
\end{align}
for all $x\in\mathbb{R}^n$. This verifies the satisfaction of Assumption \ref{assumption.steadystate}.

For the $z$-subsystem, under Assumption \ref{assumption.integralcontrol.GAS}, with the coordinate transformation \eqref{feedbackoptimization.coordinatetransformation}, we can verify properties \eqref{integralcontrol.fastsubsystem.Lya1'}--\eqref{integralcontrol.fastsubsystem.Lya3'} by defining $V_f(z,x)=V_f^o(z,x+x_*^o)$, $\underline{\alpha}_f=\underline{\alpha}_f^o$ and $\overline{\alpha}_f=\overline{\alpha}_f^o$.

For the $x$-subsystem, we consider a quadratic Lyapunov function candidate
\begin{align}
V_s(x)=\frac{1}{2}x^Tx.
\end{align}
Then, condition \eqref{assumption.slowsubsystem.Lya1} is satisfied with $\underline{\alpha}_s(r)=\overline{\alpha}_s(r)=r^2/2$ for $r\geq 0$. Moreover, under Assumptions \ref{assumption.feedbackoptimization.objectivefunction} and \ref{assumption.feedbackoptimization.Lipschitz}, direct calculation yields:
\begin{align}
&\frac{\partial V_s(x)}{\partial x}g_s(x,z)\nonumber\\
&=x^Tg_s(x,\varphi(x))+x^T(g_s(x,\varphi(x)+\tilde{z})-g_s(x,\varphi(x)))\nonumber\\
&=-x^T(\nabla\breve{\Phi}(x)-\nabla\breve{\Phi}(x_*))+x^T(g_s(x,\varphi(x)+\tilde{z})-g_s(x,\varphi(x)))\nonumber\\
&\leq-\omega|x|^2+L_{gs}|x||\tilde{z}|
\end{align}
with $\tilde{z}=z-\varphi(x)$. Thus, property \eqref{assumption.slowsubsystem.Lya2} holds with
\begin{align}
\gamma_s(r)=\frac{1}{2}\left(\frac{\mu L_{gs}}{\omega}\underline{\alpha}_f^{-1}(r)\right)^2,~~~\alpha_s(r)=2\omega\left(1-\frac{1}{\mu}\right)r,
\end{align}
and constant $\mu>1$.

Then, it can be checked that all the conditions in Assumption \ref{assumption.stability} are satisfied.

With the satisfaction of Assumptions \ref{assumption.steadystate} and \ref{assumption.stability}, Theorem \ref{theorem.tuningslow} can be readily used to find a $\rho_s^0\in\mathcal{P}$ for GAS of the feedback optimization system at the desired equilibrium.
\end{IEEEproof}

Here gives a numerical example for the feedback optimization result. With the objective function in the quadratic form, the example also verifies our result of nonlinear integral control in Section \ref{section.applications.subsection.integralcontrol}.

\begin{example}\label{example.feedbackoptimization}
Consider system
\begin{align}
\dot{z}^o&=-(z^o-x^o)^3\label{example.plant}
\end{align}
with $z^o\in\mathbb{R}$ and $x^o\in\mathbb{R}$, which is in the form of \eqref{integralcontrol.plant1}. The objective is to design an update law for $x^o$ to steer $z^o$ to the minimizing point of
\begin{align}
\Phi(z^o)=(z^o)^2-2z^o
\end{align}
by using gradient measurements. With
\begin{align}
h_f^o(z^o)=z^o,~~~\varphi^o(x^o)=x^o,
\end{align}
we consider the variable-metric gradient-flow algorithm \eqref{gradientflow} with
\begin{align}
g_s^o(x^o,z^o)=-2z^o+2.\label{example.plant}
\end{align}

In accordance with the objective function and the plant output map, set
\begin{align}
x_*^o=1,
\end{align}
and introduce a coordinate transformation:
\begin{align}
x=x^o-1,~~~z=z^o.
\end{align}
Then, we can rewrite the feedback optimization system as
\begin{align}
\dot{x}&=\rho_s^0(|-2z+2|)(-2z+2)\label{example.feedbackoptimization.plant1'}\\
\dot{z}&=-(z-x-1)^3,\label{example.feedbackoptimization.plant2'}
\end{align}
which satisfies Assumption \ref{assumption.steadystate}. The satisfaction of Assumption \ref{assumption.stability} can also be verified by defining
\begin{align}
V_s(x)=x^2,~~~V_f(z,x)=(z-x-1)^2,
\end{align}
and choosing
\begin{align}
\underline{\alpha}_s(r)&=\overline{\alpha}_s(r)=r^2,&\gamma_s(r)&=4r,&\alpha_s(r)&=2r,\\
\underline{\alpha}_f(r)&=\overline{\alpha}_f(r)=r^2,&\alpha_f(r)&=2r^2,&\lambda_f(r)&=2r^{\frac{1}{2}}.
\end{align}

Set $\underline{\rho}_f(r)\equiv 1$ and $\overline{\gamma}_s(r)=4.41r$. Then, condition \eqref{theorem.main.con3} is satisfied with
\begin{align}
\overline{\rho}_s(r)=0.99r^{\frac{3}{2}}.
\end{align}
According to \eqref{gupperbound}, we define
\begin{align}
\hat{g}_s(r)=6.2r^{\frac{1}{2}},
\end{align}
and according to \eqref{rhos0condition}, we choose
\begin{align}
\rho_s^0(r)=0.004r^2.
\end{align}

With Theorem \ref{theorem.tuningslow}, the system \eqref{example.feedbackoptimization.plant1'}--\eqref{example.feedbackoptimization.plant2'} with $(x,z)$ as the state is GAS at the equilibrium $(0,1)$, which guarantees the achievement of the feedback optimization objective.

If $\rho_s^0$ takes a positive constant value (say, $\rho_s^0(r)\equiv c>0$ for all $r\geq 0$), then the closed-loop system is reduced to
\begin{align}
\dot{x}=c(-2z+2),~~~\dot{z}=-(z-x)^3,
\end{align}
of which, the locally linearized model at the equilibrium $(x,z)=(1,1)$ is unstable for any nonzero constant $c$.
\end{example}

\subsection{Formation-Based Source Seeking}
\label{section.applications.subsection.sourceseeking}

Consider the scenario in Section \ref{section.motivation.subsection.sourceseeking}. We make the following assumption on the objective function $h$.

\begin{assumption}\label{assumption.sourceseeking.convexobjectivefunction}
There exist positive constants $\omega$ and $\vartheta$ such that
\begin{align}
(p_1-p_2)^T(\nabla h(p_1)-\nabla h(p_2))&\geq\omega|p_1-p_2|^2\label{formationsourceseekingconvexity}\\
|\nabla h(p_1)-\nabla h(p_2)|&\leq\vartheta|p_1-p_2|\label{formationsourceseekinggradientLipschitz}
\end{align}
for any $p_1,p_2\in\mathbb{R}^n$.
\end{assumption}

Property \eqref{formationsourceseekingconvexity} means that $h$ is strongly convex and admits a unique minimizer:
\begin{align}
p_*=\argmin_{p\in\mathbb{R}^n}h(p).
\end{align}
To simplify the discussion, without loss of generality, we also assume $p_*=0$.

Through an analysis based on the refined singular perturbation theorem, the following theorem shows the validity of the formation-based source-seeking algorithm in Section \ref{section.motivation.subsection.sourceseeking}.

\begin{theorem}\label{theorem.sourceseeking}
Consider a group of $N$ mobile agents modeled by \eqref{multiagent} with distributed controllers defined by \eqref{coordinatedsourceseekingcontrollaw}, \eqref{formationcontrollaw}, \eqref{formationsourceseekinglaw}, \eqref{gradientestimation1} and \eqref{gradientestimation2}. Under Assumption \ref{assumption.sourceseeking.convexobjectivefunction}, for any $p_{\epsilon}>0$, one can always find coefficients for the distributed controllers such that the closed-loop states keep bounded and the average position $p_0$ defined by \eqref{averageposition} ultimately converges to the region such that $|p_0|\leq p_{\epsilon}$.
\end{theorem}

\begin{IEEEproof}
We show that the closed-loop system can be transformed into a singularly perturbed system in the form of \eqref{plant.slow}--\eqref{plant.fast} satisfying Assumptions \ref{assumption.steadystate} and \ref{assumption.stability}. Then, Theorem \ref{theorem.mainconstant} can be used to conclude the proof.

{\em Step 1: The Closed-Loop System as a Singularly Perturbed System.} To simplify the discussion, we consider the case in which for any $i\neq j$, $a_{ij}=a_{ji}=b_{ij}=b_{ji}$ and $a_{ij}$ takes a value from $\{0,1\}$. We use an undirected graph $\mathcal{G}$ to represent the interconnection topology, with the nodes corresponding to the agents and for any $i\neq j$, the edge $(i,j)$ exists if $a_{ij}\neq 0$. It is assumed that $\mathcal{G}$ is connected.

It is a direct consequence that $\sum_{i=1}^Nv_i^f=0$, and thus the average position $p_0$ defined by \eqref{averageposition} satisfies
\begin{align}
\dot{p}_0&=-\frac{c_0}{N}\sum_{i=1}^N\sigma\left(\left(\sum_{j=1}^Nd_{j0}d_{j0}^T\right)^{-1}\delta_i\right).\label{averagepositiondynamics}
\end{align}

For $i=1,\ldots,N$, define the formation control error of agent $i$ as
\begin{align}
\tilde{p}_i=p_i-p_0-d_{i0}.
\end{align}
The boundedness of $\sigma$ guarantees the boundedness of $\tilde{p}_i$. Then, the gradient-estimation algorithm \eqref{gradientestimation1}--\eqref{gradientestimation2} can be rewritten as
\begin{align}
\dot{\delta}_i&=-(\delta_i-Nd_{i0}h(p_0+d_{i0})+w_i)-\sum_{j=1}^Na_{ij}(q_i-q_j)\label{gradientestimationalgorithm1'}\\
\dot{q}_i&=\mu\sum_{j=1}^Na_{ij}(\delta_i-\delta_j)\label{gradientestimationalgorithm2'}
\end{align}
for $i=1,\ldots,N$, with
\begin{align}
w_i=Nd_{i0}(h(p_0+d_{i0})-h(p_0+d_{i0}+\tilde{p}_i)).
\end{align}

Define
\begin{align}
\delta&=[\delta_1^T,\ldots,\delta_N^T]^T,&q&=[q_1^T,\ldots,q_N^T]^T,\\
w&=[w_1^T,\ldots,w_N^T]^T,&d_0&=[d_{10}^T,\ldots,d_{N0}^T]^T.
\end{align}

The Laplacian of the graph $\mathcal{G}$, denoted by $L$, is defined as $L=[l_{ij}]_{N\times N}$ with $l_{ii}=\sum_{j\neq i}a_{ij}$ and $l_{ij}=-a_{ij}$ for $j\neq i$. Choose $U_1\in\mathbb{R}^{N\times(N-1)}$ such that $U=[U_1,\mathbf{1}_N/\sqrt{N}]$ is a unitary matrix, and define $\hat{q}=(U_1^T\otimes I_n)q$, $\bar{L}=LU_1$ and $\bar{L}_{\otimes}=\bar{L}\otimes I_n$. The connectivity of $\mathcal{G}$ guarantees that $\bar{L}_{\otimes}^T\bar{L}_{\otimes}$ is invertible. Then, from \eqref{gradientestimationalgorithm1'}--\eqref{gradientestimationalgorithm2'}, we obtain a reduced-dimensional model of the gradient-estimation algorithm:
\begin{align}
\dot{\delta}&=-(\delta-H(p_0)+w)- \bar{L}_{\otimes}\hat{q}\label{gradientestimationalgorithm1''}\\
\dot{\hat{q}}&=\mu\bar{L}_{\otimes}^T\delta\label{gradientestimationalgorithm2''}
\end{align}
where
\begin{align}
H(p_0)=N[h(p_0+d_{10})d_{10}^T,\ldots,h(p_0+d_{N0})d_{N0}^T]^T.
\end{align}

Define
\begin{align}
\varphi(p_0)=[\mathbf{1}_N^T\otimes\delta_e^T(p_0),\hat{q}_e^T(p_0)]^T\label{sourceseekingfastequilibrium}
\end{align}
with
\begin{align}
\delta_e(p_0)&=\sum_{i=1}^Nd_{i0}h(p_0+d_{i0}),\\
\hat{q}_e(p_0)&=-\left(\bar{L}_{\otimes}^T\bar{L}_{\otimes}\right)^{-1}\bar{L}_{\otimes}^T(\delta_e(p_0)-H(p_0)).
\end{align}
Direct calculation verifies that for each fixed $p_0$, $\varphi(p_0)$ is an equilibrium of the system \eqref{gradientestimationalgorithm1''}--\eqref{gradientestimationalgorithm2''} with $w\equiv 0$.

The interconnected system composed of the average position system \eqref{averagepositiondynamics} and the gradient-esimtation system \eqref{gradientestimationalgorithm1''}--\eqref{gradientestimationalgorithm2''} is transformed into a singularly perturbed system in the form of \eqref{plant.slow}--\eqref{plant.fast} by defining
\begin{align}
x=p_0,~~~z=[\delta^T,\hat{q}^T]^T,\label{sourceseekingstatedefinition}
\end{align}
and
\begin{align}
g_s(x,z)&=-\frac{1}{N}\sum_{i=1}^N\sigma\left(\left(\sum_{j=1}^Nd_{j0}d_{j0}^T\right)^{-1}\delta_i\right),\label{sourceseekinggs}\\
g_f(z,x,w)&=\left[\begin{array}{cc}
-I_{nN} & - \bar{L}_{\otimes}\\
\mu\bar{L}_{\otimes}^T & 0
\end{array}\right]\left[\begin{array}{c}\delta\\ \hat{q}
\end{array}\right]+\left[\begin{array}{c}H(p_0)-w\\ 0
\end{array}\right].\label{sourceseekinggf}
\end{align}

{\em Step 2: Satisfaction of Assumptions \ref{assumption.steadystate} and \ref{assumption.stability}.} The satisfaction of Assumption \ref{assumption.steadystate} can be directly checked through the discussion above. Now, we verify the satisfaction of Assumption \ref{assumption.stability} by studying the stability properties of the subsystems. Define
\begin{align}
V_s(x)=\frac{1}{2}x^Tx.
\end{align}
Then, condition \eqref{assumption.slowsubsystem.Lya1} is satisfied with $\underline{\alpha}_s(s)=\overline{\alpha}_s(s)=s^2/2$. To check the satisfaction of condition \eqref{assumption.slowsubsystem.Lya2}, we rewrite $g_s$ as follows:
\begin{align}
g_s(x,z)&=-\frac{1}{N}\sum_{i=1}^N\sigma\left(\left(\sum_{j=1}^Nd_{j0}d_{j0}^T\right)^{-1}\delta_i\right)\nonumber\\
&=-\frac{1}{N}\sum_{i=1}^N\sigma\left(\left(\sum_{j=1}^Nd_{j0}d_{j0}^T\right)^{-1}(\delta_e(p_0)+\delta_i-\delta_e(p_0))\right)\nonumber\\
&=-\frac{1}{N}\sum_{i=1}^N\sigma\left(\left(\sum_{j=1}^Nd_{j0}d_{j0}^T\right)^{-1}\left(\sum_{i=1}^Nd_{i0}h(p_0+d_{i0})\right)\right.\nonumber\\
&~~~~+\left.\left(\sum_{j=1}^Nd_{j0}d_{j0}^T\right)^{-1}(\delta_i-\delta_e(p_0))\right).\label{coordinatedsourceseekinggs}
\end{align}
Moreover, with Lemma \ref{lemma.approximation}, we have
\begin{align}
\sum_{i=1}^Nd_{i0}h(p_0+d_{i0})&=\sum_{i=1}^Nd_{i0}(h(p_0)+\nabla^T h(p_0)d_{i0}+\upsilon_i(p_0))\nonumber\\
&=\sum_{i=1}^Nd_{i0}(d_{i0}^T\nabla h(p_0)+\upsilon_i(p_0))\nonumber\\
&=\left(\sum_{i=1}^Nd_{i0}d_{i0}^T\right)\nabla h(p_0)+\sum_{i=1}^Nd_{i0}\upsilon_i(p_0)
\end{align}
where $|\upsilon_i(p_0)|\leq \bar{\upsilon}|d_{i0}|^2$ with $\bar{\upsilon}$ being a positive constant. Then,
\begin{align}
g_s(x,z)&=-\frac{1}{N}\sum_{i=1}^N\sigma\left(\nabla h(p_0)+\left(\sum_{j=1}^Nd_{j0}d_{j0}^T\right)^{-1}\sum_{i=1}^Nd_{i0}\upsilon_i(p_0)\right.\nonumber\\
&~~~~+\left.\left(\sum_{j=1}^Nd_{j0}d_{j0}^T\right)^{-1}(\delta_i-\delta_e(p_0))\right)
\end{align}

Denote
\begin{align}
\Delta&=\left(\sum_{j=1}^Nd_{j0}d_{j0}^T\right)^{-1}\sum_{i=1}^Nd_{i0}\upsilon_i(p_0)+\left(\sum_{j=1}^Nd_{j0}d_{j0}^T\right)^{-1}(\delta_i-\delta_e(p_0)).
\end{align}
Consider the case in which
\begin{align}
(1-\theta)\omega|p_0|\geq|\Delta|
\end{align}
with $\theta$ being a constant satisfying $0<\theta<1$. In this case, if $p_0\neq 0$, using conditions \eqref{formationsourceseekingconvexity} and \eqref{formationsourceseekinggradientLipschitz}, we have
\begin{align}
|\nabla h(p_0)+\Delta|&\geq\frac{p_0^T}{|p_0|}(\nabla h(p_0)+\Delta)\geq\theta\omega|p_0|,\\
|\nabla h(p_0)+\Delta|&\leq(\vartheta+(1-\theta)\omega)|p_0|,
\end{align}
and thus,
\begin{align}
&p_0^T\sigma(\nabla h(p_0)+\Delta)\nonumber\\
&=p_0^T(\nabla h(p_0)+\Delta)\min\left\{1,\frac{1}{|\nabla h(p_0)+\Delta|}\right\}\nonumber\\
&\geq\left(p_0^T\nabla h(p_0)-\max_{|\Delta|\leq(1-\theta)\omega|p_0|}p_0^T\Delta\right)\min\left\{1,\frac{1}{(\vartheta+(1-\theta)\omega)|p_0|}\right\}\nonumber\\
&\geq\left(\omega|p_0|^2-(1-\theta)\omega|p_0|^2\right)\min\left\{1,\frac{1}{(\vartheta+(1-\theta)\omega)|p_0|}\right\}\nonumber\\
&=\min\left\{\theta\omega|p_0|^2,\frac{\theta\omega|p_0|}{\vartheta+(1-\theta)\omega}\right\}.
\end{align}

Define
\begin{align}
\bar{d}&=2\left(\sum_{j=1}^Nd_{j0}d_{j0}^T\right)^{-1}\sum_{i=1}^Nd_{i0}\bar{\upsilon}|d_{i0}|^2,\\
k&=2\left|\sum_{j=1}^Nd_{j0}d_{j0}^T\right|^{-1}.
\end{align}

Recall that $x=p_0$ and $V_s(x)=x^Tx/2$. The above mathematical derivation proves that
\begin{align}
&V_s(x)\geq\frac{1}{2(1-\theta)^2\omega^2}\max_{i=1,\ldots,N}\{k^2|\delta_i-\delta_e(x)|,\bar{d}^2\}\Rightarrow\frac{\partial V_s(x)}{\partial x}g_s(x,z)\leq-\alpha_s(V_s(x))
\end{align}
where
\begin{align}
\alpha_s(s)=\min\left\{2\theta\omega s,\frac{\theta\omega\sqrt{2s}}{\vartheta+(1-\theta)\omega}\right\}
\end{align}
for $s\in\mathbb{R}_+$. Recall the definitions of $z$ and $\varphi$ in \eqref{sourceseekingstatedefinition} and \eqref{sourceseekingfastequilibrium}, respectively. Condition \eqref{assumption.slowsubsystem.Lya2} can be verified as long as $V_f$ is positive definite and radially unbounded with respect to $z-\varphi(x)$.

With $z=[\delta^T,\hat{q}^T]^T$, $g_f$ defined by \eqref{sourceseekinggf} and the equilibrium map $\varphi$ defined by \eqref{sourceseekingfastequilibrium}, the dynamics of the $z$-subsystem can be rewritten as
\begin{align}
\dot{z}=A(z-\varphi(x))+Bw
\end{align}
with
\begin{align}
A=\left[\begin{array}{cc}
-I_{nN} & - \bar{L}_{\otimes}\\
\mu\bar{L}_{\otimes}^T & 0
\end{array}\right],~~~B=\left[\begin{array}{c}
-I_{nN}\\
0
\end{array}\right]
\end{align}
With $\bar{L}_{\otimes}$ being full column rank, we can prove that $A$ is Hurwitz. This means the existence of a $P=P^T>0$ such that
\begin{align}
PA+A^TP=-I_{n(2N-1)}.
\end{align}

Define
\begin{align}
V_f(z,x)=(z-\varphi(x))^TP(z-\varphi(x)).
\end{align}
Then, condition \eqref{assumption.fastsubsystem.Lya1} is satisfied with $\underline{\alpha}_f(s)=\lambda_{\min}(P)s^2$ and $\overline{\alpha}_f(s)=\lambda_{\max}(P)s^2$. Moreover,
\begin{align}
\frac{\partial V_f(z,x)}{\partial z}g_f(z,x,w)\leq-c_1V_f(z,x)+c_2|w|^2
\end{align}
where $c_1=(1-\o)/\lambda_{\max}(P)$ and $c_2=1/(4\o)$ with constant $0<\o<1$. Then, condition \eqref{assumption.fastsubsystem.Lya2} is satisfied by choosing
\begin{align}
\chi_f(s)=\frac{c_2}{c_1-c_3}s^2,~~~\alpha_f(s)= c_3 s
\end{align}
with constant $c_3$ satisfying $0<c_3<c_1$. Moreover, direct calculation gurantees the satisfaction of condition \eqref{assumption.fastsubsystem.Lya3} with $\lambda_{f1}(r)=\lambda_{f1}^cr$ and $\lambda_{f2}(r)=\lambda_{f2}^cr$ for $r\geq 0$, where $\lambda_{f1}^c$ and $\lambda_{f2}^c$ are positive constants.

With $g_s$ in the form of \eqref{coordinatedsourceseekinggs}, one may choose $\sigma$ arbitrarily small to guarantee the existence of $c_0>0$ to satisfy condition \eqref{theorem.mainconstant.condition}. With Theorem \ref{theorem.mainconstant}, it can be concluded that the closed-loop system is ISS with $\bar{d}$ and $w$ as the inputs. Moreover, from the technical derivations above, it can be recognized that $\bar{d}$ and the upper bound of $|w|$ can be rendered arbitrarily small but the corresponding gains keep unchanged by appropriately choosing the coefficients of the source-seeking algorithm. This ends the proof of Theorem \ref{theorem.sourceseeking}.
\end{IEEEproof}

\section{Conclusions}
\label{section.conclusions}

This paper has presented a generalization of singular perturbation theory, in which state-dependent perturbation functions replace perturbation coefficients, aiming to address time scales that depend on system states. This new framework has been shown to be beneficial in resolving input-to-state stabilization problems for singularly perturbed nonlinear systems. The main results are proved by using ISS methods and the nonlinear small-gain theorem. We anticipate that the new theory will be further developed to cover various practical scenarios in hierarchical, distributed, and cascade systems, and contribute to the development of emerging control architectures.

Extensions with boundary-layer subsystems possessing time-varying integral manifolds \cite{Balachandra-Sethna75,Hale80} would be beneficial for refined designs of adaptive and learning systems \cite{Ljung77,Anderson86,Solo-Kong95,Sastry-Bodson89}. The fruitful results combining averaging and singular perturbation techniques also motivate further theoretical developments with perturbation functions. Furthermore, it would be of both theoretical and practical interest to explore the generalized singular perturbation theory for systems that involve nonsmooth or hybrid dynamics, using the tools presented in \cite{Clarke-Ledyaev-Stern-Wolenski98,Goebel-Sanfelice-Teel12,Sanfelice-Teel11,Wang-Teel-Nesic12}.

\appendices

\section{Input-to-State Stability (ISS) and the Nonlinear Small-Gain Theorem}
\label{section.preliminaries}

This section reviews the basic notions and results of ISS and the nonlinear small-gain theorem. The initial development of ISS and the ISS small-gain theorem can be found in \cite{Sontag89} and \cite{Jiang-Teel-Praly94}, respectively. For tutorials, one can refer to \cite{Sontag07} and \cite{Jiang-Liu18}.

Consider a nonlinear system:
\begin{align}
\dot{x}=f(x,d)\label{generalsystem}
\end{align}
where $x\in\mathbb{R}^n$ is the state, $d\in\mathbb{R}^m$ represents the disturbance input, $f:\mathbb{R}^n\times\mathbb{R}^m\rightarrow\mathbb{R}^n$ is locally Lipschitz. It is assumed that $d$ is measurable and locally essentially bounded.

\begin{definition}
The system \eqref{generalsystem} is said to be input-to-state stable (ISS), if there exist $\beta\in\mathcal{KL}$ and $\chi\in\mathcal{K}$ such that for any initial state $x(0)$ and any input $d$,
\begin{align}
|x(t)|\leq\max\{\beta(|x(0)|,t),\chi(\|d\|_{\infty})\}
\end{align}
holds for all $t\geq 0$.
\end{definition}

\begin{theorem}
The system \eqref{generalsystem} is ISS, if it admits a continuously differentiable ISS-Lyapunov function $V:\mathbb{R}^n\rightarrow\mathbb{R}_+$, for which, there exist $\underline{\alpha},\overline{\alpha}\in\mathcal{K}_{\infty}$,  $\gamma\in\mathcal{K}$ and $\alpha\in\mathcal{PD}$ such that
\begin{align}
\underline{\alpha}(|x|)\leq V(x)\leq \overline{\alpha}(|x|)
\end{align}
for all $x\in\mathbb{R}^n$, and
\begin{align}
V(x)\geq\gamma(|d|)\Rightarrow\frac{\partial V(x)}{\partial x}f(x,d)\leq-\alpha(V(x))
\end{align}
for all $x\in\mathbb{R}^n$ and $d\in\mathbb{R}^m$.
\end{theorem}

An interconnection of two ISS subsystems with states $x_1\in\mathbb{R}^{n_1}$ and $x_2\in\mathbb{R}^{n_2}$ can be grouped together in the form of \eqref{generalsystem},
with $x=[x_1^T,x_2^T]^T$ and $n=n_1+n_2$. We assume that each subsystem admits an ISS-like Lyapunov function.

\begin{assumption}\label{assumption.generalsubsystemISS}
For $i=1,2$, there exist continuously differentiable functions $V_i:\mathbb{R}^n\rightarrow\mathbb{R}_+$ satisfying
\begin{align}
&\underline{\alpha}_i(|x_i|)\leq V_i(x)\leq \overline{\alpha}_i(|x_i|)\\
&V_i(x)\geq\max\left\{\gamma_i(V_{3-i}(x)),\chi_i(|d|)\right\}\Rightarrow\frac{\partial V_i(x)}{\partial x}f(x,d)\leq-\alpha_i(V_i(x))
\end{align}
for all $x\in\mathbb{R}^n$ and $d\in\mathbb{R}^m$, with $\underline{\alpha}_i,\overline{\alpha}_i,\gamma_i,\chi_i\in\mathcal{K}_{\infty}$ and $\alpha_i\in\mathcal{PD}$.
\end{assumption}

The Lyapunov-based ISS small-gain theorem given below plays a vital role in the Lyapunov-based small-gain analysis in this paper.

\begin{theorem}\label{theorem.Lyapunovsmallgain}
Under Assumption \ref{assumption.generalsubsystemISS}, the system \eqref{generalsystem} is ISS with $x$ as the state and $d$ as the input, if the small-gain condition
\begin{align}
\gamma_1\circ\gamma_2<\id\label{smallgaincondition}
\end{align}
is satisfied. Moreover, under the small-gain condition,
\begin{align}
V(x)=\max\left\{V_1(x),\sigma(V_2(x))\right\}
\end{align}
with $\sigma\in\mathcal{K}_{\infty}$ being continuously differentiable on $(0,\infty)$ and satisfying $\gamma_1<\sigma<\gamma_2^{-1}$, is positive definite and radially unbounded with respect to $x$, and satisfies
\begin{align}
V(x)\geq\chi(|d|)\Rightarrow\frac{\partial V(x)}{\partial x}f(x,d)\leq-\alpha(V(x))
\end{align}
wherever $V$ is differentiable, with $\chi(r)=\max\{\chi_1(r),\sigma\circ\chi_2(r)\}$ for $r\geq 0$, and $\alpha\in\mathcal{PD}$.
\end{theorem}

Slightly different from \cite{Jiang-Mareels-Wang96}, Theorem \ref{theorem.Lyapunovsmallgain} considers the case in which $V_i$ possibly depends on $x_{3-i}$. Nevertheless, the techniques in \cite{Jiang-Mareels-Wang96} are still valid to prove Theorem \ref{theorem.Lyapunovsmallgain}. 

%

\section{A Technical Lemma on First-Order Approximation of Continuously Differentiable Functions}

The following lemma is employed to develop the formation-based source-seeking algorithm using the refined singular perturbation theorem.

\begin{lemma}\label{lemma.approximation}
Consider a continuously differentiable function $h:\mathbb{R}^n\rightarrow\mathbb{R}$. Suppose that there exists a positive constant $\vartheta$ such that
\begin{align}
|\nabla h(\xi_1)-\nabla h(\xi_2)|&\leq\vartheta|\xi_1-\xi_2|\label{Lipschitzcondition}
\end{align}
for all $\xi_1,\xi_2\in\mathbb{R}^n$. Then,
\begin{align}
|h(\xi_1)-h(\xi_2)-\nabla h(\xi_1)(\xi_1-\xi_2)|\leq\vartheta|\xi_1-\xi_2|^2
\end{align}
holds for all $\xi_1,\xi_2\in\mathbb{R}^n$.
\end{lemma}

\begin{IEEEproof}
For any $\xi_1,\xi_2\in\mathbb{R}^n$, the differential mean-value theorem \cite{Rudin76} guarantees the existence of $0\leq b_i\leq 1$ for $i=1,2,\cdots,n$ such that
\begin{align}
h(\xi_1)-h(\xi_2)= \nabla h(\xi_1+(I_n-b)(\xi_1-\xi_2))(\xi_1-\xi_2)
\end{align}
with $b=\diag\{b_1,b_2,\cdots,b_n\}$. Denote $\tilde{\xi}=\xi_1-\xi_2$. Then, under condition \eqref{Lipschitzcondition}, it follows that
\begin{align}
&|h(\xi_1)-h(\xi_2)-\nabla h(\xi_1)(\xi_1-\xi_2)|\nonumber\\
&=|\nabla h(\xi_1+(I_n-b)\tilde{\xi})\tilde{\xi}-\nabla h(\xi_1)\tilde{\xi}|\nonumber\\
&\leq|\nabla h(\xi_1+(I_n-b)\tilde{\xi})-\nabla h(\xi_1)||\tilde{\xi}|\nonumber\\
&\leq\vartheta |I_n-b| |\tilde \xi|^2 \nonumber \\
&\leq \vartheta |\tilde \xi|^2.
\end{align}
This ends the proof of Lemma \ref{lemma.approximation}.
\end{IEEEproof}

\bibliographystyle{IEEEtran}
\bibliography{IEEEabrv,DFCReferences}

\end{document}